\documentclass[12pt]{article}
\usepackage{multicol}
\usepackage{amsmath, amssymb, amscd, enumerate, amsfonts, latexsym, amsxtra, 
amsthm, bm, multirow}
\usepackage{indentfirst}
\allowdisplaybreaks[3]
\newtheorem{The}{Theorem}[section]
\newtheorem{Prop}[The]{Proposition}
\newtheorem{Lem}[The]{Lemma}
\newtheorem{Cor}[The]{Corollary}
\newtheorem{Def}[The]{Definition}

\begin{document}
\centerline{\Large A basis of a certain module for the hyperalgebra} \vspace{3mm}
\centerline{\Large of $({\rm SL}_2)_r$ and some applications} \vspace{7mm}
\centerline{Yutaka Yoshii 
\footnote{ E-mail address: yutaka.yoshii.6174@vc.ibaraki.ac.jp}}  \vspace{5mm}
\centerline{College of Education,   
Ibaraki University,}
\centerline{2-1-1 Bunkyo, Mito, Ibaraki, 310-8512, Japan}
\begin{abstract}
In the hyperalgebra $\mathcal{U}_r$ of the $r$-th Frobenius kernel $({\rm SL}_2)_r$ of the 
algebraic group ${\rm SL}_2$, we construct  a basis of the $\mathcal{U}_r$-module    
 generated by a certain element which was given by the author before. 
As its applications, we also prove some results on the $\mathcal{U}_r$-modules 
and  the algebra $\mathcal{U}_r$.  
\end{abstract}

\section{Introduction}
Let $k$ be an algebraically closed 
field of characteristic $p>0$. Let $G$ be a simply connected and simple algebraic group over $k$  which is defined and split over the finite field $\mathbb{F}_p$ of 
$p$ elements. Let $G_r$ be the $r$-th Frobenius kernel of $G$. 

The representation theory of $G_r$ and $G$ is well-known. For example, as for simple modules, 
any simple $G$-module can be written in terms of some simple $G$-modules and tensor 
products by Steinberg's tensor product theorem, and any simple $G_r$-module 
can be obtained by restricting a simple $G$-module to $G_r$. Moreover, if $p$ is not too 
small,   all projective indecomposable modules (PIMs) for $G_r$ can be lifted to $G$, 
as is proved in Jantzen's paper \cite{jantzen80}. 

On the other hand, let 
$\mathcal{U}_r$ and $\mathcal{U}$ be the hyperalgebras of $G_r$ and $G$ respectively (see  
\cite[I.7.7]{jantzenbook} for definition of hyperalgebras).  
The representation theory of $G_r$ can be identified with that of $\mathcal{U}_r$, and 
the representation theory of $G$ can be identified with the locally finite one of $\mathcal{U}$ 
(see I.8.6 and II.1.20 in \cite{jantzenbook}). 
The representation theory of  $\mathcal{U}_r$ also contributes to that of $\mathcal{U}$. Since 
$\mathcal{U}_r$ is a finite-dimensional $k$-algebra, it is also 
interesting to study ring-theoretic  
information on $\mathcal{U}_r$ which contains the Jacobson radical and 
a decomposition of the unity $1 \in \mathcal{U}_r$ 
into a sum of pairwise orthogonal primitive idempotents. 
However, these are little known for  general $G$. As for the simplest case 
$G={\rm SL}_2$ (the special linear group of degree 2 over $k$), for odd $p$ and $r=1$, a generating set of the Jacobson radical of $\mathcal{U}_1$ 
is given by Wong \cite{wong83} and  a decomposition of the unity $1 \in \mathcal{U}_1$ 
as above is given by Seligman 
\cite{seligman03}.   
Recently, for any positive integer $r$, the author gives several elements in $\mathcal{U}_r$ 
(when $G={\rm SL}_2$)   denoted by $B^{({\bm \varepsilon})}({\bm a}, {\bm j})$ including  
pairwise orthogonal primitive idempotents $B^{({\bm 0})}({\bm a}, {\bm j})$ whose sum is 
the unity $1$.    

In this paper, we give some facts on the ring-theoretic 
information of the $k$-algebra $\mathcal{U}_r$ for $G={\rm SL}_2$ using the elements 
$B^{({\bm \varepsilon})}({\bm a}, {\bm j})$. 
The main results  will be given in Section 5. First we construct a basis of the 
$\mathcal{U}_r$-module $\mathcal{U}_r B^{({\bm \varepsilon})}({\bm a}, {\bm j})$ in 
Theorem 5.1. This theorem  has some applications. Indeed, it enables us to show 
that each $\mathcal{U}_r$-module $\mathcal{U}_rB^{({\bm 1})}({\bm a}, {\bm j})$ is simple 
and to construct  socle of the $\mathcal{U}_r$-module $\mathcal{U}_r$. Moreover, it also 
enables us to give a basis of the Jacobson radical of  $\mathcal{U}_r$. This result can be expected to construct a `good' generating set of its Jacobson radical for general $r$ 
which improves the main result in  \cite{wong83} when $p$ is odd and $r=1$.

\section{Preliminaries}
From now on, let $G={\rm SL}_2$. Let 
$$X= \left( \begin{array}{cc}
         0 & 1 \\
         0 & 0 \end{array} \right),\ \ \ 
Y= \left( \begin{array}{cc}
         0 & 0 \\
         1 & 0 \end{array} \right),\ \ \ 
H= \left( \begin{array}{cc}
         1 & 0 \\
         0 & {-1} \end{array} \right)$$
be the standard basis in the simple complex Lie algebra $\mathfrak{g}_{\mathbb{C}}=\mathfrak{sl}_2(\mathbb{C})$. Let  
$\mathcal{U}_{\mathbb{Z}}$ be the subring of the universal enveloping algebra $\mathcal{U}_{\mathbb{C}}$ of  $\mathfrak{g}_{\mathbb{C}}$  
generated by all $X^{(m)}=X^m/m!$ and  $Y^{(m)}=Y^m/m!$ with $m \in \mathbb{Z}_{\geq 0}$. 
For  $z=H$ or $-H$ in $\mathcal{U}_{\mathbb{Z}}$, set 
$${z+c \choose m }= \dfrac{(z+c)(z+c-1) \cdots (z+c-m+1)}{m!}$$
for $c \in \mathbb{Z}$ and $m \in \mathbb{Z}_{\geq 0}$, which also lies in 
$\mathcal{U}_{\mathbb{Z}}$. Then the elements 
$$Y^{(m)} {H \choose n} X^{(m')}$$
with $m, m', n \in \mathbb{Z}_{\geq 0}$ form a $\mathbb{Z}$-basis of $\mathcal{U}_{\mathbb{Z}}$. The $k$-algebra 
$\mathcal{U}_{\mathbb{Z}} \otimes_{\mathbb{Z}} k$ is denoted by $\mathcal{U}$ and is 
identified with  
the hyperalgebra of $G$, which is equipped with a structure of Hopf algebra over $k$. 
We  use the same notation for the images in $\mathcal{U}$ 
of the elements 
in $\mathcal{U}_{\mathbb{Z}}$.

To carry out calculations in $\mathcal{U}$, we use the following well-known 
equalities:  
$$\displaystyle{X^{(m)} Y^{(n)} = \sum_{i=0}^{{\rm min}(m,n)} Y^{(n-i)} 
{H-m-n+2i \choose i} X^{(m-i)}},$$ 
$$\displaystyle{Y^{(m)} X^{(n)} = \sum_{i=0}^{{\rm min}(m,n)} X^{(n-i)} 
{-H-m-n+2i \choose i} Y^{(m-i)}},$$
$$\displaystyle{{H +s \choose m} X^{(n)}= X^{(n)} {H+s +2n \choose m}}, \ \ \   \displaystyle{{H +s \choose m} Y^{(n)}= Y^{(n)} {H+s -2n \choose m}},$$
$$\displaystyle{X^{(m)} X^{(n)} = {m+n \choose n} X^{(m+n)}}, \ \ \   
\displaystyle{Y^{(m)} Y^{(n)} = {m+n \choose n} Y^{(m+n)}},$$
$$\displaystyle{{H \choose m} {H \choose n}= \sum_{i=0}^{{\rm min}(m,n) }
\dfrac{(m+n-i)!}{(m-i)! (n-i)! i!} {H \choose m+n-i}},$$
$$\displaystyle{{H+s+t \choose m} = 
\sum_{i=0}^{m} {s \choose m-i} {H+t \choose i}}$$
for $m,n \in \mathbb{Z}_{\geq 0}$ and  
$s,t \in \mathbb{Z}$.

Let 
${\rm Fr}: \mathcal{U} \rightarrow \mathcal{U}$ be the $k$-algebra endomorphism 
 defined by 
$${\rm Fr} \left(X^{(m)}\right) = 
\left\{ \begin{array}{cl}
         {X^{(m/p)}} & {\mbox{if $p \mid m$,}} \\
         0 & {\mbox{if $p \nmid m$}} \end{array} \right. 
\ \ \ \  \mbox{and  \ \ \ \ }
{\rm Fr} \left(Y^{(m)}\right) = 
\left\{ \begin{array}{cl}
         {Y^{(m/p)}} & {\mbox{if $p \mid m$,}} \\
         0 & {\mbox{if $p \nmid m$}} \end{array} \right. .$$
Then we also have 
$${\rm Fr} \bigg( {H \choose m} \bigg) = 
\left\{ \begin{array}{cl}
         {{H \choose m/p}} & {\mbox{if $p \mid m$,}} \\
         0 & {\mbox{if $p \nmid m$}} \end{array} \right. .$$

The elements 
$Y^{(m)} {H \choose n} X^{(m')}$ 
with $m, m', n \in \mathbb{Z}_{\geq 0}$ form a $k$-basis of $\mathcal{U}$. 
We say that a nonzero  element 
$z \in \mathcal{U}$ has degree $d$ if it is a $k$-linear combination of  the elements 
$Y^{(m)} {H \choose n} X^{(m')}$ 
with $m, m', n \in \mathbb{Z}_{\geq 0}$ and $m'-m =d$.
For a positive integer 
$r \in \mathbb{Z}_{> 0}$, let $\mathcal{U}_r$ be the subalgebra of $\mathcal{U}$ generated by 
$X^{(n)}$ and $Y^{(n)}$ with $0 \leq n  \leq p^r-1$. 
This is a finite-dimensional $k$-algebra of dimension $p^{3r}$ which has 
$Y^{(m)} {H \choose n} X^{(m')}$ 
with $0 \leq m, m', n \leq p^r-1$ as a basis, 
and it can be identified with the hyperalgebra  of the $r$-th Frobenius kernel 
$G_r = {\rm Ker} ({\rm Fr}^r)$ of $G$. 
Let $\mathcal{U}^0$ be 
the subalgebra of $\mathcal{U}$ generated by ${H \choose n}$ with 
$n \in \mathbb{Z}_{\geq 0}$ and set 
$\mathcal{U}_r^0=\mathcal{U}^0 \cap \mathcal{U}_r$.  
The subalgebra  $\mathcal{U}_r^0$ is generated by ${H \choose n}$ with 
$0 \leq n  \leq p^r-1$.

Let  ${\rm Fr}': \mathcal{U} \rightarrow \mathcal{U}$ be the $k$-linear map defined 
by 
$$Y^{(m)} {H \choose n} X^{(m')} \mapsto 
Y^{(mp)} {H \choose np} X^{(m'p)}.$$
This map is used to construct the elements 
$B^{({\bm \varepsilon})}({\bm a},{\bm j})$ defined in the next section.    

Let $\mathcal{A}$ be the $k$-subalgebra of $\mathcal{U}$ which is generated by 
$\mathcal{U}^0$ and all $Y^{(n)}X^{(n)}$ with $n \in \mathbb{Z}_{\geq 0}$. 
The subalgebra $\mathcal{A}$ is commutative and consists of the zero 
element and all elements of 
degree 0 in $\mathcal{U}$. 
For a positive integer $r \in \mathbb{Z}_{>0}$, set 
$\mathcal{A}_r=\mathcal{A} \cap \mathcal{U}_r$. 
For any $n \in \mathbb{Z}_{\geq 0}$, the elements  $X^{(pn)}$ and 
$Y^{(pn)}$ commute with all elements in $\mathcal{A}_1$. 
For details, see \cite[\S 2]{yoshii17}.

Throughout this paper, all modules for an associative 
 $k$-algebra are assumed to be finite-dimensional left 
modules and $r$ denotes a fixed positive integer unless otherwise stated. 
For a finite-dimensional (associative) $k$-algebra $R$, let ${\rm rad}R$ be the 
largest nilpotent two-sided ideal of $R$, which is called the Jacobson radical of $R$. 
Let $M$ be an $R$-module. Then the $R$-submodule $({\rm rad} R) M$ is denoted by 
${\rm rad}_R M$ and called the radical of $M$. The quotient 
$R$-module $M/ {\rm rad}_R M$ is called the head of $M$. On the other hand,   
the $R$-submodule of $M$ consisting of the elements annihilated by ${\rm rad} R$ 
is denoted by ${\rm soc}_R M$ and called the  socle of $M$. 
This is  the largest semisimple $R$-submodule of $M$. For details, for example, see 
\cite[\S 1]{alperinbook}.

\section{Construction of the elements 
$B^{({\bm \varepsilon})}({\bm a},{\bm j})$ in $\mathcal{U}_r$} 
In this section we describe a construction of some elements in $\mathcal{A}_r$,  
which are denoted by $B^{({\bm \varepsilon})}({\bm a},{\bm j})$ 
and include pairwise orthogonal  primitive idempotents in $\mathcal{U}_r$ whose sum is the unity $1 \in \mathcal{U}_r$. For details, see 
 \cite[\S 4 and 5]{yoshii17} and \cite[\S 3]{yoshii18}. 

For $a \in \mathbb{Z}$, set 
$$\mu_a = {H-a-1 \choose p-1} = \sum_{i=0}^{p-1}
{-a-1 \choose p-1-i} {H \choose i} \in \mathcal{U}_1^0.$$
This is a   $\mathcal{U}_1^0$-weight vector of weight $a$ in the $\mathcal{U}_1^0$-module  
$\mathcal{U}_1^0$: $H \mu_a = a \mu_a$.
Moreover,  we have  $\mu_a = \mu_b$ if and only if $a \equiv b\ ({\rm mod}\ p)$, and 
all $\mu_a$ with $a \in \{ 0,1, \dots, p-1\}$ are pairwise orthogonal 
primitive idempotents in $\mathcal{U}_1^0$ whose sum is the unity $1 \in \mathcal{U}_1^0$.

Suppose for a moment  that $p$ is odd. 
Set $\mathcal{S} = \{ 0, 1, \dots , (p-1)/2 \} \subset \mathbb{Z}$. 
We denote by $\overline{\mathcal{S}}$  the image of the subset 
$\mathcal{S} \subset \mathbb{Z}$ under the natural map 
$\mathbb{Z} \rightarrow \mathbb{F}_p$. 
For $\varepsilon \in \mathbb{F}_2 = \{ 0,1 \}$ and $j \in \mathcal{S}$,  
 we define  polynomials $\psi(x), \psi_j^{(\varepsilon)}(x) \in \mathbb{F}_p [x]$  as  
$$\psi(x)=\prod_{i \in \mathbb{F}_p} \left(x- i^2\right) = x 
\prod_{i \in \overline{\mathcal{S}} \backslash \{ 0\}} 
\left(x- i^2\right)^2, $$
$$\psi_0^{(0)}(x)= \psi_0^{(1)}(x)= \prod_{i \in \mathbb{F}_p^{\times}} \left(x-i^2\right) = 
\prod_{i \in \overline{\mathcal{S}} \backslash \{ 0\}} \left(x-i^2\right)^2,$$
$$\psi_{s}^{(0)}(x) = 2x\left(x+s^2\right) \prod_{i \in \mathbb{F}_p^{\times} 
\backslash \{ s, -s \} } \left(x-i^2\right)
 = 2x\left(x+s^2\right) \prod_{i \in \overline{\mathcal{S}} \backslash \{0, s \} } 
\left(x-i^2\right)^2,  
$$
and 
$$\psi_s^{(1)}(x) =x (x-s^2)
\prod_{i \in \overline{\mathcal{S}} \backslash \{ 0,s\}} 
\left(x- i^2\right)^2 $$
if $s \in \mathcal{S} \backslash \{ 0\}$ 
($i$ and $s$ in the right-hand sides denote their images under the natural map 
$\mathbb{Z} \rightarrow \mathbb{F}_p$). 
Set $\mathcal{P}_{\mathbb{Z}} = \mathbb{Z} \times \mathcal{S}$
and 
$$B^{(\varepsilon)}(a, j) = \psi_j^{(\varepsilon)} 
\left(\mu_a YX +\left(\dfrac{a+1}{2}\right)^2\right) \cdot \mu_a$$
for $\varepsilon \in \mathbb{F}_2$ and $(a,j) \in \mathcal{P}_{\mathbb{Z}}$.   
Since $\psi_0^{(0)}(x)=\psi_0^{(1)}(x)$, we have $B^{(0)}(a,0)=B^{(1)}(a,0)$ for any 
$a \in \mathbb{Z}$. 

In turn, suppose that  $p=2$. Then we consider the set 
$$\mathcal{P}_{\mathbb{Z}}=\left\{ \left. \left(2i,\dfrac{1}{2}\right), (1+2i,0), (1+2i,1) 
\ \ \right|\   
i \in \mathbb{Z} \right\} \subset 
\mathbb{Z} \times \mathbb{Q}$$ instead of $\mathcal{P}_{\mathbb{Z}}
=\mathbb{Z} \times \mathcal{S}$ when $p$ is odd, and define 
$$B^{(0)}\left(2i,\dfrac{1}{2}\right)= \mu_0,\ \ \ B^{(1)}\left(2i,\dfrac{1}{2}\right)= 
\mu_0YX= \mu_0 XY, $$
$$B^{(0)}(1+2i,0)=B^{(1)}(1+2i,0)= \mu_1YX = \mu_1 XY + \mu_1,$$
$$B^{(0)}(1+2i,1)=  B^{(1)}(1+2i,1)= \mu_1YX + \mu_1=\mu_1XY$$
for any $i \in \mathbb{Z}$.

Let $p$ be an arbitrary prime number again. 
For $\varepsilon \in \mathbb{F}_2$ and 
$(a_1,j_1), (a_2,j_2) \in \mathcal{P}_{\mathbb{Z}}$, we have  
$$B^{(\varepsilon)}(a_1,j_1) = B^{(\varepsilon)}(a_2,j_2) 
\Longleftrightarrow \mbox{$a_1 \equiv a_2$ $({\rm mod}\ p)$ and $j_1=j_2$.}$$  
We define a subset $\mathcal{P}$ of $\mathcal{P}_{\mathbb{Z}}$ as 
$\mathcal{P} = \{ (a,j) \in \mathcal{P}_{\mathbb{Z}}\ |\ 
0 \leq a \leq p-1\}$, and hence 
$\mathcal{P}= \{ 0,1, \dots, p-1 \} \times \mathcal{S}$ if $p$ is odd and 
$\mathcal{P}=\{(0,1/2), (1,0), (1,1) \}$ if $p=2$. 
The elements $B^{(0)}(a,j)$ with $(a,j) \in \mathcal{P}$ are 
 pairwise orthogonal primitive idempotents in $\mathcal{U}_1$ whose sum is   the unity $1 \in \mathcal{U}_1$.    

For an integer $n \in \mathbb{Z}$, we denote by $n\ {\rm {\bf mod}}\ p$  
a unique integer $\widehat{n}$ with 
$\widehat{n} \equiv n\ ({\rm mod}\ p)$ and $0 \leq \widehat{n} \leq p-1$. 

We classify  pairs $(a,j) \in \mathcal{P}_{\mathbb{Z}}$  
under the following four conditions: \\ \\
(A) $\widehat{a}$ is even and $(p-\widehat{a}+1)/2 \leq j \leq (p-1)/2$,\\
(B) $\widehat{a}$ is even and $0 \leq j \leq (p-\widehat{a}-1)/2$,\\
(C) $\widehat{a}$ is odd and $0 \leq j \leq (\widehat{a}-1)/2$,\\
(D) $\widehat{a}$ is odd and $(\widehat{a}+1)/2 \leq j \leq (p-1)/2$,\\ \\
where $\widehat{a}= a\ {\rm {\bf mod}}\ p$. 
Note that if $p=2$, the pairs $(2i,1/2)$, $(1+2i,0)$, and $(1+2i,1)$ in 
$\mathcal{P}_{\mathbb{Z}}$  
for $i \in \mathbb{Z}$ satisfy (B), (C), and (D) 
respectively.  Apart from them,  we also consider the following condition for 
$(a,j) \in \mathcal{P}_{\mathbb{Z}}$: \\ \\
(E) $j=0$ if $p$ is odd or $a \equiv 1\ ({\rm mod}\ 2)$ if $p=2$. \\

\begin{Def}
Let $\varepsilon \in \mathbb{F}_2$ and  $(a,j) \in \mathcal{P}_{\mathbb{Z}}$, and 
set $\widehat{a}= a\ {\rm {\bf mod}}\ p$. Then 
define  nonnegative integers 
$n^{(\varepsilon)}(a,j)$ and 
$\widetilde{n}^{(\varepsilon)}(a,j)$ every condition of $(a,j)$ 
from {\rm (A)} to {\rm (D)} as follows: 
$$\begin{array}{|c||c|c|c|c|} \hline 
& & & & \\[-3mm]
(a,j) & n^{(0)}(a,j) & n^{(1)}(a,j) & \widetilde{n}^{(0)}(a,j) & \widetilde{n}^{(1)}(a,j) \\[1mm] \hline \hline 
& & & & \\[-3mm] 
{\rm (A)} & \dfrac{p-\widehat{a}-1}{2}+j   & \dfrac{3p-\widehat{a}-1}{2}-j & 
\dfrac{-p+\widehat{a}-1}{2}+j & \dfrac{p+\widehat{a}-1}{2}-j 
\\[3mm] \hline 
& & & & \\[-3mm]
{\rm (B)} & \dfrac{p-\widehat{a}-1}{2}-j & \dfrac{p-\widehat{a}-1}{2}+j & 
\dfrac{p+\widehat{a}-1}{2}-j & \dfrac{p+\widehat{a}-1}{2}+j 
\\[3mm] \hline 
& & & & \\[-3mm]
{\rm (C)} & \dfrac{2p-\widehat{a}-1}{2}-j & \dfrac{2p-\widehat{a}-1}{2}+j & 
\dfrac{\widehat{a}-1}{2}-j & \dfrac{\widehat{a}-1}{2}+j 
\\[3mm] \hline 
& & & & \\[-3mm]
{\rm (D)} & j-\dfrac{\widehat{a}+1}{2} & \dfrac{2p-\widehat{a}-1}{2}-j & 
\dfrac{\widehat{a}-1}{2}+j & \dfrac{2p+\widehat{a}-1}{2}-j 
\\[3mm] \hline 
\end{array}$$ 
\end{Def}
\ \\

\noindent {\bf Remark.} For $(a,j) \in \mathcal{P}_{\mathbb{Z}}$ and 
$\varepsilon \in \mathbb{F}_2$, we easily see the following. \\ 

\noindent {\rm (a)} $\widetilde{n}^{(\varepsilon)}(a,j)=n^{(\varepsilon)}(-a,j)$. \\ 

\noindent {\rm (b)} $0 \leq n^{(0)}(a,j) \leq n^{(1)}(a,j) \leq p-1$ and 
$$n^{(0)}(a,j)=n^{(1)}(a,j) \Longleftrightarrow 
\mbox{$(a,j)$ satisfies (E)}. $$

\noindent {\rm (c)}
$n^{(0)}(a,j)+ \widetilde{n}^{(1)}(a,j)=n^{(1)}(a,j)+ \widetilde{n}^{(0)}(a,j)=p-1$. \\

\begin{Lem}
Let $(a,j) \in \mathcal{P}_{\mathbb{Z}}$ and $\varepsilon \in \mathbb{F}_2$. Then the element $B^{(\varepsilon)}(a,j)$ 
 can be written as
$$B^{(\varepsilon)}(a,j)= \mu_a \sum_{m=n^{(\varepsilon)}(a,j)}^{p-1} 
c^{(\varepsilon)}_m(a,j) Y^m X^m = 
\mu_a \sum_{m=\widetilde{n}^{(\varepsilon)}(a,j)}^{p-1} 
\widetilde{c}^{(\varepsilon)}_m(a,j) X^m Y^m$$
for some $c^{(\varepsilon)}_m(a,j), \widetilde{c}^{(\varepsilon)}_m(a,j) \in \mathbb{F}_p$ with 
$c^{(\varepsilon)}_{n^{(\varepsilon)}(a,j)}(a,j) \neq 0$ and 
$\widetilde{c}^{(\varepsilon)}_{\widetilde{n}^{(\varepsilon)}(a,j)}(a,j) \neq 0$. 
\end{Lem}
\ 

Recall that $r$ is a fixed positive integer. For  
$a \in \mathbb{Z}$, set 
$$\mu_a^{(r)} = {H-a-1 \choose p^r-1} \in \mathcal{U}^0_{r}.$$
Clearly we have $\mu_a^{(1)}=\mu_a$. \\

\begin{Prop}
For $a \in \mathbb{Z}$, the following hold. \\ 

\noindent {\rm (i)} $\mu_a^{(r)}$ is a   $\mathcal{U}_r^0$-weight vector of weight $a$ in the $\mathcal{U}_r^0$-module  
$\mathcal{U}_r^0$:  
${H \choose n} \mu^{(r)}_a = {a \choose n} \mu^{(r)}_a$ 
for any $n \in \{0,1,\dots, p^r-1 \}$. \\

\noindent {\rm (ii)} For $b \in \mathbb{Z}$, we have 
$$\mu^{(r)}_a = \mu^{(r)}_b \Longleftrightarrow \mbox{$a \equiv b\ ({\rm mod}\ p^r)$}.$$
\ 

\noindent {\rm (iii)} The elements 
 $\mu^{(r)}_a$ with $a \in \{ 0,1, \dots, p^r-1\}$ are pairwise orthogonal 
primitive idempotents in $\mathcal{U}_r^0$ whose sum is  the unity $1 \in \mathcal{U}_r^0$. \\

\noindent {\rm (iv)} We have 
$$\mu^{(r)}_a X^{(n)}= X^{(n)}\mu^{(r)}_{a-2n}\ \ \ \mbox{and}\ \ \ 
\mu^{(r)}_a Y^{(n)}= Y^{(n)}\mu^{(r)}_{a+2n} $$
for any $n \in \mathbb{Z}_{\geq 0}$. \\

\noindent {\rm (v)} Suppose that $r \geq 2$. 
If $a=a'+p^i a''$ with $0 \leq a' \leq p^i-1$, $a'' \in \mathbb{Z}$, and 
$1 \leq i \leq r-1$, we have 
$$\mu^{(r)}_a= \mu_{a'}^{(i)} {\rm Fr}'^i \left( \mu_{a''}^{(r-i)}\right).$$
\end{Prop}
\ 
 
\noindent For details, see \cite[\S4]{gros-kaneda15}.

If $(a,j) \in \mathcal{P}_{\mathbb{Z}}$  satisfies 
{\rm (A)} or {\rm (C)}, then define an integer $s(a,j)$ as 
$$s(a,j)=\dfrac{p-(a \ {\rm {\bf mod}}\ p)+1}{2}$$
 if $p$ is odd and 
$a\ {\rm {\bf mod}}\ p$ is even,  
$$s(a,j)=\dfrac{p-(a\ {\rm {\bf mod}}\ p)}{2}$$ if both $p$ and $a\ {\rm {\bf mod}}\ p$ are  odd, 
and  $s(a,j)=1$ if $p=2$. 

For $\varepsilon \in \mathbb{F}_2$ and $(a,j) \in \mathcal{P}_{\mathbb{Z}}$,  we write 
$$B^{(\varepsilon)}(a,j)= \mu_a \sum_{m=n^{(\varepsilon)}(a,j)}^{p-1} 
c^{(\varepsilon)}_m(a,j) Y^m X^m$$
as in  Lemma 3.2. Then    
we define $Z^{(\varepsilon)} \left(z; (a,j) \right)$ for $z \in \mathcal{U}$ as 
$$Z^{(\varepsilon)} \left(z; (a,j) \right)
= \mu_a \sum_{m=n^{(\varepsilon)}(a,j)}^{p-1} 
c^{(\varepsilon)}_m(a,j) Y^m X^{m-s(a,j)} {\rm Fr'}(z)X^{s(a,j)}$$
if $(a,j)$ satisfies (A) or (C), and
$$Z^{(\varepsilon)} \left(z; (a,j) \right)
={\rm Fr'}(z) B^{(\varepsilon)}(a,j)\ \left( =B^{(\varepsilon)}(a,j){\rm Fr'}(z)\right)$$
if $(a,j)$ satisfies (B) or (D). \\

\begin{Prop} 
Let $(a,j) \in \mathcal{P}_{\mathbb{Z}}$ and $\varepsilon \in \mathbb{F}_2$. 
Then the following hold. \\ \\
{\rm (i)} The map $Z^{(\varepsilon)} \left(-; (a,j) \right) : \mathcal{U} \rightarrow \mathcal{U}, \ 
z \mapsto Z^{(\varepsilon)} \left(z; (a,j) \right)$ is  $k$-linear and injective. \\ \\
{\rm (ii)} For  $z \in \mathcal{U}$,  there is  
an element $z' \in \mathcal{U}$ which is independent of 
$\varepsilon$ such that $$Z^{(\varepsilon)} \left(z; (a,j) \right)
={\rm Fr'}(z') B^{(\varepsilon)}(a,j)=B^{(\varepsilon)}(a,j){\rm Fr'}(z').$$
Then we also have 
$$z=0 \Longleftrightarrow {\rm Fr}'(z')=0 \Longleftrightarrow z'=0$$
and 
$$z \in \mathcal{A} \Longleftrightarrow Z^{(\varepsilon)} \left(z; (a,j) \right) \in \mathcal{A}
\Longleftrightarrow z' \in \mathcal{A}.$$
\ 

\noindent {\rm (iii)} Let $u$ be an element of the $k$-subalgebra of $\mathcal{U}$ generated by all $X^{(p^i)}$ and $Y^{(p^i)}$ 
with $i \in \mathbb{Z}_{> 0}$ (and the unity $1 \in \mathcal{U}$). 
Then we have 
$u Z^{(\varepsilon)} \left(z; (a,j) \right)= 
Z^{(\varepsilon)} \left( {\rm Fr}(u)z; (a,j) \right)$. \\ \\
{\rm (iv)} For any $z_1,z_2 \in \mathcal{U}$, we have 
$$Z^{(0)} \left(z_1; (a,j) \right) Z^{(\varepsilon)} \left(z_2; (a,j) \right)
=Z^{(\varepsilon)} \left(z_1; (a,j) \right) Z^{(0)} \left(z_2; (a,j) \right)
= Z^{(\varepsilon)} \left(z_1z_2; (a,j) \right).$$ \\ 
{\rm (v)} Suppose that $(a,j)$ satisfies the condition {\rm (E)}. 
Then for  any 
$z \in \mathcal{U}$, we have 
$$Z^{(0)} \left(z; (a,j) \right)=Z^{(1)} \left(z; (a,j) \right).$$
\end{Prop} 
\ 

\noindent {\bf Remark.} 
(a) \ In (ii), the equivalences associated to the condition 
$z=0$ follow from the injectivity of the maps 
$Z^{(\varepsilon)} \left( -;(a,j)\right)$ and ${\rm Fr}'$ and from the fact that the multiplication 
map $\mathcal{U}_1 \otimes_k {\rm Fr}'(\mathcal{U}) \rightarrow \mathcal{U}$ is a 
$k$-linear isomorphism (see the remark of \cite[Proposition 2.3]{yoshii17}). 

(b) \ (iv) is a generalization of the fact that 
the map $Z^{(0)} \left(-; (a,j) \right)$ is multiplicative 
(see \cite[Proposition 5.4 (ii)]{yoshii17}). The equalities are obtained 
by letting $B^{(\varepsilon)}(a,j)$ act on the equality 
$$Z^{(0)} \left(z_1; (a,j) \right) Z^{(0)} \left(z_2; (a,j) \right)
= Z^{(0)} \left(z_1z_2; (a,j) \right)$$
from the left and the right respectively. \\

Consider an $r$-tuple 
$\left((a_i, j_i)\right)_{i=0}^{r-1} =\left( (a_0,j_0),\dots, (a_{r-1},j_{r-1}) \right)
 \in \mathcal{P}_{\mathbb{Z}}^r$ 
of  pairs 
$(a_i,j_i) \in \mathcal{P}_{\mathbb{Z}}$ $(0 \leq i \leq r-1)$. For convenience we shall write it 
 as 
$$((a_0,  \dots, a_{r-1}),(j_0,  \dots, j_{r-1}))$$
or $({\bm a},{\bm j})$ with ${\bm a}=(a_0,  \dots, a_{r-1})$ and 
${\bm j}=(j_0,  \dots, j_{r-1})$.

For ${\bm \varepsilon}=(\varepsilon_0, \dots, \varepsilon_{r-1}) \in \mathbb{F}_2^r$ and  
$({\bm a},{\bm j})= \left((a_i, j_i)\right)_{i=0}^{r-1} \in \mathcal{P}_{\mathbb{Z}}^r$,  
we define an element 
$B^{({\bm \varepsilon})}({\bm a},{\bm j}) = 
B^{(\varepsilon_0, \dots, \varepsilon_{r-1})}
\left((a_0,  \dots, a_{r-1}),(j_0,  \dots, j_{r-1})\right)
\in \mathcal{U}$ inductively 
as $B^{(\varepsilon_0)}(a_0,j_0) $ if $r=1$, and 
$$
Z^{(\varepsilon_0)} \left( 
B^{(\varepsilon_1, \dots, \varepsilon_{r-1})}\left((a_1,  \dots, a_{r-1}),(j_1,  \dots, j_{r-1})\right); (a_0,j_0) \right)
$$
if $r \geq 2$.  \\

\begin{Prop} Let $({\bm a}, {\bm j}) = 
\left( (a_i,j_i)\right)_{i=0}^{r-1}\in \mathcal{P}_{\mathbb{Z}}$ and 
${\bm \varepsilon} =(\varepsilon_0, \dots, \varepsilon_{r-1}) \in \mathbb{F}_2^r$. 
Then the following hold. \\ 

\noindent {\rm (i)}  $B^{(\bm{\varepsilon} )}({\bm a},{\bm j})$ lies in $\mathcal{A}_r$ and is 
a $\mathcal{U}_r^0$-weight vector of 
 $\mathcal{U}_r^0$-weight $\sum_{i=0}^{r-1}p^i b_i $, where 
$$b_i = 
\left\{ \begin{array}{ll} 
a_i\ {\rm {\bf mod}}\ p-p & {\mbox{if $(a_i,j_i)$ 
satisfies {\rm (A)} or {\rm (C)},}} \\
a_i\ {\rm {\bf mod}}\ p   & {\mbox{if $(a_i,j_i)$ 
satisfies {\rm (B)} or {\rm (D)}}}
\end{array} \right.. $$

\noindent {\rm (ii)} Let ${\bm \rho} = (\rho_0, \dots, \rho_{r-1}) \in \mathbb{F}_2^r$. Then 
$B^{({\bm \varepsilon})}({\bm a},{\bm j})=B^{({\bm \rho})}({\bm a},{\bm j})$ if and only if 
$\varepsilon_i = \rho_i$ whenever $(a_i,j_i)$ does not satisfy the condition {\rm (E)}. \\  

\noindent {\rm (iii)} Set ${\bm 0}=(0, \dots, 0) \in \mathbb{F}_2^{r}$. 
The elements $B^{({\bm 0})}({\bm a}, {\bm j})$ with $({\bm a}, {\bm j}) \in \mathcal{P}^{r}$ 
are pairwise orthogonal primitive 
idempotents in $\mathcal{U}_r$ whose sum is the unity $1 \in \mathcal{U}_r$. 
\end{Prop}
\ 

\noindent {\bf Remark.} (iii) implies that the sum 
$\sum_{({\bm a}, {\bm j}) \in\mathcal{P}^r} \mathcal{U}_r  
B^{({\bm 0})}({\bm a}, {\bm j})$  
gives a direct sum decomposition of the $\mathcal{U}_r$-module $\mathcal{U}_r$ into PIMs  
(see \cite[ch. 1. Theorem 4.7]{nagao-tsushimabook}).

\section{Preparation for main results}

In this paper, we need a lot of notation for main results. So we collect them here. \\

\begin{Def}
Let $(a,j) \in \mathcal{P}_{\mathbb{Z}}$ and 
$({\bm a}, {\bm j})= \left( (a_i,j_i)\right)_{i=0}^{r-1} \in \mathcal{P}_{\mathbb{Z}}^r$. \\ 

\noindent {\rm (1)} For $i \in \mathbb{Z}$ and $n \in \{0,1, \dots, p-1 \}$, define 
$\gamma_{i}(a,j)$, $\widetilde{\gamma}_{i}(a,j)$, $\beta_{n}(a,j)$, and 
$\widetilde{\beta}_{n}(a,j)$ in $\mathbb{F}_p$ as follows: 
$$\gamma_{i}(a,j)= j^2- \left( \dfrac{a+1}{2} \right)^2 -i (i+a+1) \ \left( =
 j^2- \left( \dfrac{a+1}{2} +i \right)^2 \right),$$
$$\widetilde{\gamma}_{i}(a,j)= \gamma_{i}(-a,j),$$
$$\beta_{n}(a,j)= \prod_{i=0}^{n-1}\gamma_{i}(a,j),$$
$$\widetilde{\beta}_{n}(a,j) = \beta_{n}(-a,j) \ \left( 
=\prod_{i=0}^{n-1}\widetilde{\gamma}_{i}(a,j) \right).$$  
Here if $p=2$, $\gamma_{i}(a,j)$ is defined by regarding the right-hand side (which is 
an integer in this situation) as the image under the natural map
$\mathbb{Z} \rightarrow \mathbb{F}_2$. If $p$ is odd, 
$\gamma_{i}(a,j)$ is defined by regarding the integers $i$, $j$, and $a+1$ in the 
right-hand side as the images 
under the natural map $\mathbb{Z} \rightarrow \mathbb{F}_p$. \\

\noindent {\rm (2)} We regard $\mathbb{F}_2^r$ as an additive group induced by 
the addition in $\mathbb{F}_2$. Define two elements ${\bm 0}$ and ${\bm 1}$ in 
$\mathbb{F}_2^r$ as 
$${\bm 0}= (0,\dots, 0),$$
$${\bm 1}= (1,\dots, 1).$$ 
Moreover, for ${\bm \varepsilon}=(\varepsilon_0, \dots, \varepsilon_{r-1})$, 
$\widetilde{\bm \varepsilon}=(\widetilde{\varepsilon}_0, \dots, \widetilde{\varepsilon}_{r-1})
 \in \mathbb{F}_2^r$,  
define   
${\bm \varepsilon} \leq \widetilde{\bm \varepsilon}$ 
 if $\varepsilon_i \leq \widetilde{\varepsilon}_i$ for each $i$,   
regarding $\varepsilon_i$ and $\widetilde{\varepsilon}_i$ as the corresponding integers 
(i.e. $0$ or $1$ in $\mathbb{Z}$). This 
gives a partial order in $\mathbb{F}_2^{r}$. \\ 

\noindent {\rm (3)} Two subsets $\mathcal{X}_r({\bm a},{\bm j})$ and 
$\mathcal{Y}_r({\bm a}, {\bm j})$ of $\mathbb{F}_2^r$ are defined as 
$$\mathcal{X}_r({\bm a}, {\bm j}) = 
\{ (\varepsilon_0, \dots, \varepsilon_{r-1})\in \mathbb{F}_2^r\ |\ 
\varepsilon_i=0\ \mbox{whenever $(a_i,j_i)$ satisfies {\rm (E)}}\}, $$
$$\mathcal{Y}_r({\bm a}, {\bm j}) = 
\{ (\varepsilon_0, \dots, \varepsilon_{r-1})\in \mathbb{F}_2^r\ |\ 
\varepsilon_i=1\ \mbox{whenever  $(a_i,j_i)$ satisfies {\rm (E)}}\}.$$
\ 

\noindent {\rm (4)} Under the partial order defined in {\rm (2)}, let 
${\bm \sigma}({\bm a},{\bm j})$ be a unique maximal element in 
 $\mathcal{X}_r({\bm a},{\bm j})$ and ${\bm \tau}({\bm a},{\bm j})$ a unique 
minimal one in $\mathcal{Y}_r({\bm a},{\bm j})$. In other words, 
${\bm \sigma}({\bm a},{\bm j})=(\sigma_0, \dots, \sigma_{r-1})$ and 
${\bm \tau}({\bm a},{\bm j})=(\tau_0, \dots, \tau_{r-1})$, where 
$$\sigma_i = \left\{ 
\begin{array}{ll}
0 & \mbox{if $(a_i,j_i)$ satisfies {\rm (E)},} \\
1 & \mbox{otherwise}
\end{array}
\right.$$
and 
$$\tau_i = \left\{ 
\begin{array}{ll}
1 & \mbox{if $(a_i,j_i)$ satisfies {\rm (E)},} \\
0 & \mbox{otherwise}
\end{array}
\right..$$ 
From now on we adopt such notation for the entries of ${\bm \sigma}({\bm a},{\bm j})$ and 
${\bm \tau}({\bm a},{\bm j})$ unless otherwise stated. \\ 

\noindent {\rm (5)} For ${\bm \varepsilon} \in \mathcal{X}_r({\bm a},{\bm j})$ 
and $\widetilde{\bm \varepsilon} \in \mathcal{Y}_r({\bm a},{\bm j})$, define two subsets 
$\Theta_r\left( ({\bm a}, {\bm j}), {\bm \varepsilon}\right)$ and 
$\widehat{\Theta}_r\left( ({\bm a}, {\bm j}), \widetilde{\bm \varepsilon}\right)$ of 
$\mathbb{F}_2^r \times \mathbb{Z}^r$ as 
$$\Theta_r\left( ({\bm a}, {\bm j}), {\bm \varepsilon}\right)=
\left\{ \left( {\bm \theta}, {\bm t}({\bm \theta})\right)\ \left|\  
\begin{array}{l}
\mbox{${\bm \varepsilon} \leq {\bm \theta} \in \mathcal{X}_r({\bm a},{\bm j})$ and} \\
\mbox{$-\widetilde{n}^{(\theta_i+1)}(a_i,j_i) \leq t_i(\theta_i) \leq n^{(\theta_i+1)}(a_i,j_i)$ 
for each $i$}
\end{array}
\right. \right\},$$
$$ \widehat{\Theta}_r\left( ({\bm a}, {\bm j}), \widetilde{\bm \varepsilon}\right)=
\left\{ \left( {\bm \theta}, {\bm t}({\bm \theta})\right)\ \left|\  
\begin{array}{l}
\mbox{$\widetilde{\bm \varepsilon} \leq {\bm \theta} \in \mathcal{Y}_r({\bm a},{\bm j})$ 
and} \\
\mbox{$-\widetilde{n}^{(\theta_i+1)}(a_i,j_i) \leq t_i(\theta_i) \leq n^{(\theta_i+1)}(a_i,j_i)$ 
for each $i$}
\end{array}
\right. \right\},$$
where ${\bm \theta}= (\theta_0, \dots, \theta_{r-1}) \in \mathbb{F}_2^r$ and 
${\bm t}({\bm \theta})= \left( t_0(\theta_0), \dots, t_{r-1}(\theta_{r-1})\right) 
\in \mathbb{Z}^r$.  
From now on we adopt such notation for the entries of 
${\bm \theta}$ and ${\bm t}({\bm \theta})$ with respect to an element 
$\left( {\bm \theta}, {\bm t}({\bm \theta})\right)$ in 
$\Theta_r\left( ({\bm a}, {\bm j}), {\bm \varepsilon}\right)$ or 
$\widehat{\Theta}_r\left( ({\bm a}, {\bm j}), \widetilde{\bm \varepsilon}\right)$ 
unless otherwise stated. 
\\ 

\noindent {\rm (6)} For $i \in \mathbb{Z}_{\geq 0}$ and 
$t \in \mathbb{Z}$, define an element 
$u^{(i, t)}$ in $\mathcal{U}$ as 
$$u^{(i, t)}=\left\{ \begin{array}{ll} 
 {X^{(p^i)t}} & {\mbox{if $t \geq 0$,}} \\ 
 {\left(Y^{(p^i)}\right)^{-t}} & {\mbox{if $t < 0$}}
\end{array} \right. .$$
Moreover, for   
${\bm \varepsilon} =(\varepsilon_0, \dots, \varepsilon_{r-1}) 
\in \mathbb{F}_2^r$ and ${\bm t}= (t_0 ,\dots, t_{r-1}) \in \mathbb{Z}^r$, 
define an element 
$B^{({\bm \varepsilon})}\left( ({\bm a}, {\bm j}); {\bm t} \right)$ in $\mathcal{U}_r$ as 
$$B^{({\bm \varepsilon})}\left( ({\bm a}, {\bm j}); {\bm t} \right)
=u^{(0,t_0)} u^{(1, t_1)} \cdots u^{({r-1}, t_{r-1})}
B^{({\bm \varepsilon})}({\bm a}, {\bm j}).$$

\noindent {\rm (7)} For ${\bm \varepsilon} \in \mathcal{X}_r({\bm a},{\bm j})$ 
and $\widetilde{\bm \varepsilon} \in \mathcal{Y}_r({\bm a},{\bm j})$, define two subsets 
$\mathcal{B}_r\left( ({\bm a}, {\bm j}), {\bm \varepsilon}\right)$ and 
$\widehat{\mathcal{B}}_r\left( ({\bm a}, {\bm j}), \widetilde{\bm \varepsilon}\right)$ of 
$\mathcal{U}_r$ as 
$$\mathcal{B}_r\left( ({\bm a}, {\bm j}), {\bm \varepsilon}\right)
= \left\{ \left. B^{({\bm \theta})}\left(({\bm a}, {\bm j}); {\bm t}({\bm \theta})\right)
\ \right|\ \left({\bm \theta}, {\bm t}({\bm \theta})\right) \in 
\Theta_{r}\left( ({\bm a}, {\bm j}), {\bm \varepsilon}\right)\right\},$$
$$\widehat{\mathcal{B}}_r\left( ({\bm a}, {\bm j}), \widetilde{\bm \varepsilon}\right)
= \left\{ B^{({\bm \theta})}\left(({\bm a}, {\bm j}); {\bm t}({\bm \theta})\right)
\ \left|\ \left({\bm \theta}, {\bm t}({\bm \theta})\right) \in 
\widehat{\Theta}_{r}\left( ({\bm a}, {\bm j}), \widetilde{\bm \varepsilon}\right) \right. \right\}.$$
\end{Def} 
\ 

\noindent {\bf Remark.} (a) Clearly we have $\beta_0(a,j)= \widetilde{\beta}_0(a,j)=1$ 
by definition. 
For $i \in \mathbb{Z}$ and $s \in \mathbb{Z}_{\geq 0}$, we 
have 
$$\gamma_{i+s}(a,j) = \gamma_{i}(a+2s,j),$$
$$\widetilde{\gamma}_{i+s}(a,j)= \widetilde{\gamma}_{i}(a-2s,j)$$
by definition. Moreover, if $0 \leq i \leq p-1$, we have 
$$\gamma_{i}(a,j)=0 \Longleftrightarrow \mbox{$i$ is equal to 
$n^{(0)}(a,j)$ or $n^{(1)}(a,j)$,} $$
$$\widetilde{\gamma}_{i}(a,j)=0 \Longleftrightarrow \mbox{$i$ is equal to 
$\widetilde{n}^{(0)}(a,j)$ or $\widetilde{n}^{(1)}(a,j)$.} $$
\ 

\noindent {\rm (b)} Each of $\mathcal{X}_r({\bm a},{\bm j})$ and 
$\mathcal{Y}_r({\bm a},{\bm j})$ is used to remove duplicates from the 
elements $B^{({\bm \varepsilon})}({\bm a}, {\bm j})$ with 
${\bm \varepsilon} \in \mathbb{F}_2^r$. Indeed, for each 
${\bm \varepsilon} \in \mathbb{F}_2^r$, there exist unique 
elements ${\bm \theta} \in \mathcal{X}_r({\bm a},{\bm j})$ and 
${\bm \rho} \in \mathcal{Y}_r({\bm a},{\bm j})$ such that 
$$B^{({\bm \varepsilon})}({\bm a},{\bm j})= 
B^{({\bm \theta})}({\bm a},{\bm j})=B^{({\bm \rho})}({\bm a},{\bm j})$$
by Proposition 3.5 (ii). We also 
note that $\mathcal{X}_r({\bm a},{\bm j})$ and $\mathcal{Y}_r({\bm a},{\bm j})$ 
can be written as 
$$\mathcal{X}_r({\bm a},{\bm j}) = 
\{ (\theta_0, \dots, \theta_{r-1}) \in \mathbb{F}_2^r\ |\ \theta_i \in \{ 0,\sigma_i \},\ 
0 \leq i \leq r-1\},$$ 
$$\mathcal{Y}_r({\bm a},{\bm j}) = 
\{ (\theta_0, \dots, \theta_{r-1}) \in \mathbb{F}_2^r\ |\ \theta_i \in \{ \tau_i,1 \},\ 
0 \leq i \leq r-1\}$$
under the notation in $(4)$.  
\\ 

\noindent {\rm (c)}  We have to keep in mind that the symbol 
${\bm \tau}= {\bm \tau}({\bm a},{\bm j})$ in \cite[\S 3]{yoshii18} denotes 
${\bm \sigma}({\bm a},{\bm j})$ (not ${\bm \tau}({\bm a},{\bm j})$) here. \\

\noindent {\rm (d)}  Clearly  
${\bm \sigma}({\bm a},{\bm j})+{\bm \tau}({\bm a},{\bm j})= {\bm 1}$. Moreover, 
for ${\bm \varepsilon} \in \mathcal{X}_r({\bm a},{\bm j})$ and 
$\widetilde{\bm \varepsilon} \in \mathcal{Y}_r({\bm a},{\bm j})$, we have 
$$B^{({\bm \varepsilon})}( {\bm a},{\bm j}) = 
B^{(\widetilde{\bm \varepsilon})}( {\bm a},{\bm j})  
\Longleftrightarrow \widetilde{\bm \varepsilon}= 
{\bm \varepsilon}+{\bm \tau}({\bm a},{\bm j})$$
by Proposition 3.5 (ii). In particular, we have 
$$B^{({\bm 0})}({\bm a},{\bm j})=B^{({\bm \tau}({\bm a},{\bm j}))}({\bm a},{\bm j}), $$  
$$B^{({\bm 1})}({\bm a},{\bm j})=B^{({\bm \sigma}({\bm a},{\bm j}))}({\bm a},{\bm j}), $$ 
and 
$$\mathcal{B}_r \left( ({\bm a},{\bm j}), {\bm \varepsilon}\right)
=\widehat{\mathcal{B}}_r \left( ({\bm a},{\bm j}), {\bm \varepsilon}+ 
{\bm \tau}({\bm a},{\bm j}) \right)$$
for ${\bm \varepsilon} \in \mathcal{X}_r({\bm a},{\bm j})$.
\\ 

\noindent {\rm (e)} Let 
${\bm \varepsilon} =(\varepsilon_0, \dots, \varepsilon_{r-1}) \in  \mathbb{F}_2^r$ 
and $\widetilde{{\bm \varepsilon}} = (\widetilde{\varepsilon}_0, \dots, \widetilde{\varepsilon}_{r-1}) \in \mathcal{X}_r({\bm a},{\bm j})$. Then the product 
$B^{({\bm \varepsilon })} ({\bm a},{\bm j})B^{(\widetilde{\bm \varepsilon })} ({\bm a},{\bm j})$ 
is equal to zero if there is an integer $s$ with $0 \leq s \leq r-1$ such that  
$\varepsilon_s = \widetilde{\varepsilon}_s =1$, and to 
$B^{({\bm \varepsilon }+\widetilde{\bm \varepsilon })} ({\bm a},{\bm j})$ otherwise. This fact is 
a generalization of \cite[Lemma 3.9]{yoshii18} and can be proved in a similar way. \\

\noindent {\rm (f)} In (6), clearly we have $u^{(i,t)}=0$ if 
$t \leq -p$ or $t \geq p$. It follows that  
$B^{({\bm \varepsilon})}\left( ({\bm a}, {\bm j}); {\bm t} \right)=0$ if there exists 
an integer $i \in \{ 0, \dots, r-1\}$ such that  $t_i \leq -p$ or $t_i \geq p$. \\

\noindent {\rm (g)} For ${\bm \varepsilon}, {\bm \rho} \in 
\mathcal{X}_r({\bm a},{\bm j})$, we have 
$$\mathcal{U}_r B^{({\bm \varepsilon })} ({\bm a},{\bm j}) \subseteq 
\mathcal{U}_r B^{({\bm \rho})} ({\bm a},{\bm j}) \Longleftrightarrow 
{\bm \varepsilon } \geq {\bm \rho},$$
$$\mathcal{U}_r B^{({\bm \varepsilon })} ({\bm a},{\bm j}) = 
\mathcal{U}_r B^{({\bm \rho})} ({\bm a},{\bm j}) \Longleftrightarrow 
{\bm \varepsilon } = {\bm \rho}.$$
This fact easily follows from (e) and implies that these equivalences also hold for 
${\bm \varepsilon}, {\bm \rho} \in 
\mathcal{Y}_r({\bm a},{\bm j})$.  
\\

\noindent {\bf Example.} Consider the case of $r=1$. 
If $(a,j)$ satisfies {\rm (E)}, we  have 
$$\mathcal{X}_1(a,j)=\{0\},\ \ \ \ \mathcal{Y}_1(a,j)=\{1\},$$
$$\Theta_1\left( (a,j), 0\right)= \left\{  \left( 0, t(0)\right)\ \left| \ 
-\widetilde{n}^{(1)}(a,j) \leq t(0) \leq n^{(1)}(a,j)\right. \right\},$$
$$\widehat{\Theta}_1\left( (a,j), 1\right)= \left\{  \left( 1, t(1)\right)\ \left| \ 
-\widetilde{n}^{(0)}(a,j) \leq t(1) \leq n^{(0)}(a,j)\right. \right\},$$
\begin{align*}
\mathcal{B}_1\left( (a,j), 0\right) &= \left\{  B^{(0)}\left( (a,j); t(0)\right)\ \left| \ 
-\widetilde{n}^{(1)}(a,j) \leq t(0) \leq n^{(1)}(a,j)\right. \right\} \\
&= \left\{  B^{(1)}\left( (a,j); t(1)\right)\ \left| \ 
-\widetilde{n}^{(0)}(a,j) \leq t(1) \leq n^{(0)}(a,j) \right. \right\} \\
&= \widehat{\mathcal{B}}_1\left( (a,j), 1\right). 
\end{align*} 
If $(a,j)$ does not satisfy {\rm (E)}, we have 
$$\mathcal{X}_1(a,j)=\mathcal{Y}_1(a,j)=\mathbb{F}_2,$$
\begin{align*}
\Theta_1\left( (a,j), 0\right) &= \left\{  \left( 0, t(0)\right), \left( 1, t(1)\right)\ \left| \ 
\begin{array}{l}
{-\widetilde{n}^{(1)}(a,j) \leq t(0) \leq n^{(1)}(a,j),} \\
{-\widetilde{n}^{(0)}(a,j) \leq t(1) \leq n^{(0)}(a,j)}
\end{array} \right. \right\} \\
&= \widehat{\Theta}_1\left( (a,j), 0\right),
\end{align*}
$$
\Theta_1\left( (a,j), 1\right) = \left\{  \left( 1, t(1)\right)\ \left| \ 
-\widetilde{n}^{(0)}(a,j) \leq t(1) \leq n^{(0)}(a,j)\right. \right\} \\
= \widehat{\Theta}_1\left( (a,j), 1\right),
$$
\begin{align*}
\mathcal{B}_1\left( (a,j), 0\right) &= \left\{  B^{(0)}\left((a,j);t(0)\right), 
B^{(1)}\left((a,j); t(1)\right)\ \left| \ 
\begin{array}{l}
{-\widetilde{n}^{(1)}(a,j) \leq t(0) \leq n^{(1)}(a,j),} \\
{-\widetilde{n}^{(0)}(a,j) \leq t(1) \leq n^{(0)}(a,j)}
\end{array} \right. \right\} \\
&= \widehat{\mathcal{B}}_1\left( (a,j), 0\right),   
\end{align*}
$$
\mathcal{B}_1\left( (a,j), 1\right) = \left\{  
B^{(1)}\left((a,j); t(1)\right)\ \left| \ 
-\widetilde{n}^{(0)}(a,j) \leq t(1) \leq n^{(0)}(a,j)\right. \right\} \\
= \widehat{\mathcal{B}}_1\left( (a,j), 1\right).   
$$
\

To describe the main results, we make some preparation. 

The proposition below describes commutation formulas between 
$B^{(\varepsilon)}(a,j)$ and $X$ or $Y$. We have to notice that 
these formulas depend on  whether $p$ is odd or even. 
\\

\begin{Prop}
Let $(a,j) \in \mathcal{P}_{\mathbb{Z}}$ and $\varepsilon \in \mathbb{F}_2$. 
Then the following hold. \\ \\
{\rm (i)} Suppose that $p$ is odd. 
Let $s,t$ be  integers with $1 \leq s \leq p-1$ and $1 \leq t \leq p-1$. Then 
$$X^s B^{(\varepsilon)}(a,j) = B^{(\varepsilon)}(a+2s,j) X^s\ \ \ \mbox{and}\ \ \ \ 
Y^t B^{(\varepsilon)}(a,j) = B^{(\varepsilon)}(a-2t,j) Y^t. $$
{\rm (ii)} Suppose that $p=2$. Then 
$$
XB^{(\varepsilon)}\left(0, \dfrac{1}{2}\right) =
B^{(\varepsilon)}\left(0, \dfrac{1}{2}\right) X, \ 
XB^{(\varepsilon)}(1,0)= B^{(\varepsilon)}(1,1)X, \   
XB^{(\varepsilon)}(1,1)= B^{(\varepsilon)}(1,0)X,$$
$$YB^{(\varepsilon)}\left(0, \dfrac{1}{2}\right)=
B^{(\varepsilon)}\left(0, \dfrac{1}{2}\right) Y,\ 
YB^{(\varepsilon)}(1,0)= B^{(\varepsilon)}(1,1)Y,\    
YB^{(\varepsilon)}(1,1)= B^{(\varepsilon)}(1,0)Y.$$ 
\end{Prop}
\ 

\noindent {\itshape Proof.} It is easy to check (ii), so we shall check 
only (i). Suppose that $p$ is odd. 
It is enough to 
check the equalities only for $s=1$ and $t=1$ in (i). Note that 
\begin{eqnarray*}
X  \left(\mu_a YX +\left(\dfrac{a+1}{2}\right)^2\right)  \mu_a  
&=& X \mu_{a} YX + \left(\dfrac{a+1}{2}\right)^2 X \mu_{a}  \\
&=& \mu_{a+2} (YX+H)X + \left(\dfrac{a+1}{2}\right)^2 \mu_{a+2} X \\
&=& \mu_{a+2} YX^2+(a+2) \mu_{a+2}X + \left(\dfrac{a+1}{2}\right)^2 \mu_{a+2} X \\
&=& \left(\mu_{a+2} YX +\left(\dfrac{a+3}{2}\right)^2\right)  \mu_{a+2} X 
\end{eqnarray*}
by Proposition 3.3 (i) and (iv). Thus we have  
\begin{eqnarray*}
X B^{(\varepsilon)}(a,j) 
&=& X \psi_j^{(\varepsilon)} 
\left(\mu_a YX +\left(\dfrac{a+1}{2}\right)^2\right)  \mu_a \\
&=& \psi_j^{(\varepsilon)} 
\left(\mu_{a+2} YX +\left(\dfrac{a+3}{2}\right)^2\right)  \mu_{a+2} X \\
&=& B^{(\varepsilon)}(a+2,j) X.
\end{eqnarray*}
Similarly, we obtain $YB^{(\varepsilon)}(a,j) =B^{(\varepsilon)}(a-2,j)Y $. $\square$ \\

\begin{Prop}
For $(a,j) \in \mathcal{P}_{\mathbb{Z}}$, the following hold. \\ 

\noindent {\rm (i)} Let $s,t$ be integers with $0 \leq s \leq n^{(0)}(a,j)$ and 
$0 \leq t \leq \widetilde{n}^{(0)}(a,j)$. Then
$$Y^s X^s B^{(0)}(a,j)= \beta_s(a,j) B^{(0)}(a,j) +
4j^2 \sum_{i=0}^{s-1} \dfrac{\beta_s(a,j)}{\gamma_i(a,j)} B^{(1)}(a,j), $$
$$X^t Y^t B^{(0)}(a,j)= \widetilde{\beta}_t(a,j) B^{(0)}(a,j) +
4j^2 \sum_{i=0}^{t-1} \dfrac{\widetilde{\beta}_t(a,j)}{\widetilde{\gamma}_i(a,j)} B^{(1)}(a,j).$$
\ 

\noindent {\rm (ii)} Let $s,t$ be integers with $0 \leq s \leq n^{(0)}(a,j)$ and 
$0 \leq t \leq \widetilde{n}^{(0)}(a,j)$. Then
$$Y^s X^s B^{(1)}(a,j)= \beta_s(a,j) B^{(1)}(a,j),$$ 
$$X^t Y^t B^{(1)}(a,j)= \widetilde{\beta}_t(a,j) B^{(1)}(a,j).$$
\ 

\noindent {\rm (iii)} Let $s,t$ be integers with $n^{(0)}(a,j) < s \leq n^{(1)}(a,j)$ and 
$\widetilde{n}^{(0)}(a,j)< t \leq \widetilde{n}^{(1)}(a,j)$ (these occur only if 
$(a,j)$ does not satisfy {\rm (E)}). Then
$$Y^s X^s B^{(0)}(a,j) = 4j^2 
\left( \prod_{i=0,\ i \neq {n}^{(0)}(a,j)}^{s-1}\gamma_i(a,j) \right) B^{(1)}(a,j),$$
$$X^t Y^t B^{(0)}(a,j) = 4j^2 
\left( \prod_{i=0,\ i \neq \widetilde{n}^{(0)}(a,j)}^{t-1}
\widetilde{\gamma}_i(a,j) \right) B^{(1)}(a,j).$$ 
\end{Prop}

\noindent {\itshape Proof.} It is easy when $p=2$, so we may assume that $p$ is odd.  We shall check the equalities only for $s$, 
since the results for  $t$ are similar. It is clear for $s=0$, so we may assume $s>0$. Note that $$\mu_aY^s X^s= 
\left( \prod_{i=0}^{s-m-1}\left(\mu_a YX - (m+i)(m+i+a+1)\right) \right) 
\cdot \mu_a Y^mX^m$$
for $0 \leq m \leq s-1$,  
$$YXB^{(1)}(a,j)=\gamma_0(a,j)B^{(1)}(a,j),$$ 
and 
$$YXB^{(0)}(a,j)=\gamma_0(a,j)B^{(0)}(a,j)+4j^2B^{(1)}(a,j)$$ in $\mathcal{U}_1$ 
(see \cite[Lemma 4.2]{yoshii17} and \cite[Proposition 3.7]{yoshii18}).  Then the first equality in (i) for $s=1$ is clear and 
the first one in (ii) is easy to check: 
\begin{align*}
Y^s X^s B^{(1)}(a,j) &=  \mu_aY^s X^s B^{(1)}(a,j)
= \left( \prod_{i=0}^{s-1}\left(\mu_a YX - i(i+a+1)\right) 
\right) B^{(1)}(a,j) \\
&= \left( \prod_{i=0}^{s-1}\left(j^2 - \left( \dfrac{a+1}{2}\right)^2 
- i(i+a+1)\right) \right) B^{(1)}(a,j) \\
&= \beta_s(a,j) B^{(1)}(a,j).
\end{align*}
Now we shall check the first equality in (i) for $s \geq 2$. 
By induction on $s$, we have 
\begin{align*}
\lefteqn{Y^s X^s B^{(0)}(a,j) =  \mu_a Y^s X^s B^{(0)}(a,j)} \\
& =  \left(\mu_a YX-(s-1)(s+a)\right) \cdot \mu_a Y^{s-1}X^{s-1} B^{(0)}(a,j) \\
& =  \left(\mu_a YX-(s-1)(s+a)\right) 
\left( \beta_{s-1}(a,j) B^{(0)}(a,j) + 
4j^2 \sum_{i=0}^{s-2} 
\dfrac{\beta_{s-1}(a,j)}{\gamma_i(a,j)} B^{(1)}(a,j )\right) \\
& =  \left(j^2 - \left( \dfrac{a+1}{2}\right)^2 - (s-1)(s+a)\right) \beta_{s-1}(a,j)B^{(0)}(a,j)
         +4j^2\beta_{s-1}(a,j)B^{(1)}(a,j) \\
& \ \ \   + 4j^2 \left(j^2 - \left( \dfrac{a+1}{2}\right)^2 - (s-1)(s+a)\right)
\sum_{i=0}^{s-2} \dfrac{\beta_{s-1}(a,j)}{\gamma_i(a,j)} 
B^{(1)}(a,j ) \\
& =  \beta_s(a,j) B^{(0)}(a,j) +
4j^2 \left( \sum_{i=0}^{s-2} 
\dfrac{\beta_{s-1}(a,j)}{\gamma_i(a,j)}\gamma_{s-1}(a,j) +
\beta_{s-1}(a,j)\right) B^{(1)}(a,j ) \\
& = \beta_s(a,j) B^{(0)}(a,j) +
4j^2\sum_{i=0}^{s-1} \dfrac{\beta_{s-1}(a,j) \gamma_{s-1}(a,j)}{\gamma_i(a,j)} B^{(1)}(a,j ) \\
& = \beta_s(a,j) B^{(0)}(a,j) +
4j^2\sum_{i=0}^{s-1} \dfrac{\beta_{s}(a,j)}{\gamma_i(a,j)} B^{(1)}(a,j ), 
\end{align*}
as required.  

If $n^{(0)}(a,j) < s \leq n^{(1)}(a,j)$, we have 
\begin{align*}
\lefteqn{Y^s X^s B^{(0)}(a,j)=\mu_a Y^s X^s B^{(0)}(a,j)} \\
& =  \left( \prod_{i=0}^{s-n^{(0)}(a,j)-1} \left(\mu_a YX - \left(n^{(0)}(a,j)+i \right)
\left(n^{(0)}(a,j)+i +a+1\right)\right) \right) \\
& \ \ \  \times \mu_a Y^{n^{(0)}(a,j)}X^{n^{(0)}(a,j)}B^{(0)}(a,j) \\
& =  \left( \prod_{i=0}^{s-n^{(0)}(a,j)-1} \left(\mu_a YX - \left(n^{(0)}(a,j)+i 
\right) \left(n^{(0)}(a,j)+i +a+1\right)\right) \right)\\
& \ \ \  \times \left( \beta_{n^{(0)}(a,j)}(a,j)B^{(0)}(a,j) +
4j^2\sum_{l=0}^{n^{(0)}(a,j)-1} \dfrac{\beta_{n^{(0)}(a,j)}(a,j)}{\gamma_l(a,j)}B^{(1)}(a,j) \right)
\end{align*}
by (i). Since 
$$\left(\mu_a YX - n^{(0)}(a,j) 
\left(n^{(0)}(a,j) +a+1\right)\right)B^{(0)}(a,j) = 4j^2 B^{(1)}(a,j),$$
$$\left(\mu_a YX - n^{(0)}(a,j) 
\left(n^{(0)}(a,j) +a+1\right)\right)B^{(1)}(a,j) =0,$$
and 
$$\left(\mu_a YX - \left(n^{(0)}(a,j)+i \right)
\left(n^{(0)}(a,j)+i +a+1\right)\right)B^{(1)}(a,j) = 
\gamma_{n^{(0)}(a,j)+i}(a,j)B^{(1)}(a,j)$$
for $1 \leq i \leq s-n^{(0)}(a,j)-1$, 
we obtain 
\begin{align*}
Y^s X^s B^{(0)}(a,j)
& =  4j^2 \beta_{n^{(0)}(a,j)}(a,j) 
\left( \prod_{i=1}^{s-n^{(0)}(a,j)-1}\gamma_{n^{(0)}(a,j)+i}(a,j) \right)
B^{(1)}(a,j) \\
& =   4j^2 \left( \prod_{i=0}^{n^{(0)}(a,j)-1}\gamma_{i}(a,j) \right)
\left( \prod_{i=n^{(0)}(a,j)+1}^{s-1}\gamma_{i}(a,j) \right)
B^{(1)}(a,j) \\
& =   4j^2 
\left( \prod_{i=0,\ i \neq n^{(0)}(a,j)}^{s-1}\gamma_{i}(a,j) \right)
B^{(1)}(a,j), 
\end{align*}
and (iii) follows. $\square$ \\

\noindent {\bf Remark.} We also note that  $X^s B^{(1)}(a,j)= 0$ if 
$n^{(0)}(a,j) < s \leq p-1$ and $Y^t B^{(1)}(a,j)=0 $ if 
$\widetilde{n}^{(0)}(a,j)< t \leq p-1$, and that $X^s B^{(0)}(a,j)= 0$ if 
$n^{(1)}(a,j) < s \leq p-1$ and $Y^t B^{(0)}(a,j)=0 $ if 
$\widetilde{n}^{(1)}(a,j)< t \leq p-1$ by Lemma 3.2 and Remark (c) of 
Definition 3.1. Since 
$\gamma_{n^{(\varepsilon)}(a,j)}(a,j)=
\widetilde{\gamma}_{\widetilde{n}^{(\varepsilon)}(a,j)}(a,j)=0$ 
for $\varepsilon \in \mathbb{F}_2$, 
each equality in (ii)  holds even if 
$n^{(0)}(a,j) < s \leq p-1$ or $\widetilde{n}^{(0)}(a,j)< t \leq p-1$ and 
each one in (iii) holds even if 
$n^{(1)}(a,j) < s \leq p-1$ or $\widetilde{n}^{(1)}(a,j)< t \leq p-1$. 
Indeed, in these cases both sides of the equalities are zero 
(in $\mathcal{U}$).   \\ 

\begin{Prop} 
Let $n \in \mathbb{Z}_{> 0}$ and     
$(a,j) \in \mathcal{P}_{\mathbb{Z}}$.    
Let $s,t$ be integers with  $0 \leq s \leq n^{(0)}(a,j)$ and 
$0 \leq t \leq \widetilde{n}^{(0)}(a,j)$. Then
$$Y^{(pn)} X^s B^{(1)}(a,j)=X^s Y^{(pn)}  B^{(1)}(a,j),$$ 
$$X^{(pn)} Y^t B^{(1)}(a,j)=Y^t X^{(pn)}  B^{(1)}(a,j).$$ 
\end{Prop}

\noindent {\itshape Proof.} It is easy when $p=2$, so we may assume that $p$ is odd. 
We shall check  the equality only for $s$. We may assume that 
$s \neq 0$. 

Suppose that $1 \leq s \leq n^{(0)}(a,j)$. Then we have 
$$Y^{(pn)} X^s B^{(1)}(a,j)=s! \sum_{i=0}^{s}X^{(s-i)} 
{-H-pn-s+2i \choose i} Y^{(pn-i)} B^{(1)}(a,j).$$
If $1 \leq i \leq s$, we have 
$Y^{(pn-i)} B^{(1)}(a,j)=Y^{(p(n-1))} Y^{(p-i)} B^{(1)}(a,j)=0$ 
by Lemma 3.2, since 
$$p-i  \geq p-s > p-1-{n}^{(0)}(a,j) \geq p-1-{n}^{(1)}(a,j).$$
Therefore, we obtain $Y^{(pn)} X^s B^{(1)}(a,j)=X^s Y^{(pn)}  
B^{(1)}(a,j)$, and the proposition follows. $\square$ \\

\noindent {\bf Remark.} Each  equality in the proposition 
holds even if 
$n^{(0)}(a,j) < s \leq p-1$ or  
$\widetilde{n}^{(0)}(a,j) < t \leq p-1$.  Indeed, in these cases the both sides 
are zero.  \\

Recall the elements $u^{(i, t)}$ and 
$B^{({\bm \varepsilon})} \left( ({\bm a}, {\bm j}) ; (t_0, \dots, t_{r-1})\right)$ in 
$\mathcal{U}$ defined in Definition 4.1 (6). \\

\begin{Prop} 
Let $({\bm a}, {\bm j})= \left( (a_i,j_i)\right)_{i=0}^{r-1} \in \mathcal{P}_{\mathbb{Z}}^r$,  
${\bm \varepsilon}=(\varepsilon_0 ,\dots, \varepsilon_{r-1}) \in \mathbb{F}_2^r$, and 
$(t_0, \dots, t_{r-1}) \in \mathbb{Z}^r$. Then 
$$B^{({\bm \varepsilon})}\big( ({\bm a}, {\bm j}); (t_0, \dots, t_{r-1}) \big) \neq 0 
\Longleftrightarrow \mbox{$-\widetilde{n}^{(\varepsilon_i+1)}(a_i,j_i) \leq t_i \leq n^{(\varepsilon_i+1)}(a_i,j_i)$ for each  $i$.}
$$
\end{Prop}

\noindent {\itshape Proof.} It is clear for $r=1$ by Lemma 3.2 
and Remark (c) of 
Definition 3.1. Assume that $r \geq 2$. 
By Proposition 3.4 (ii) note that 
$$Z^{(\varepsilon_0)}
\left( B^{({\bm \varepsilon}')}\left( ({\bm a}', {\bm j}') ; 
(t_1, \dots, t_{r-1}) \right) ; (a_0,j_0)\right)
= B^{(\varepsilon_0)}(a_0,j_0) {\rm Fr'}(z)$$
for some $z \in \mathcal{U}$, where 
$({\bm a}', {\bm j}') = \left( (a_i,j_i)\right)_{i=1}^{r-1} 
\in \mathcal{P}_{\mathbb{Z}}^{r-1}$ and 
${\bm \varepsilon}'=(\varepsilon_1, \dots, \varepsilon_{r-1}) \in 
\mathbb{F}_2^{r-1}$. Then by Proposition 3.4 (iii) we have   
\begin{align*}
B^{({\bm \varepsilon})}\left( ({\bm a}, {\bm j}); (t_0, \dots, t_{r-1}) \right)
&= u^{(0,t_0)} u^{(1,t_1)} \cdots u^{({r-1},t_{r-1})}
B^{({\bm \varepsilon})}({\bm a}, {\bm j}) \\
&= u^{(0,t_0)} u^{(1,t_1)} \cdots u^{({r-1},t_{r-1})} Z^{(\varepsilon_0)}
\left( B^{({\bm \varepsilon}')}({\bm a}', {\bm j}') ; (a_0,j_0)\right) \\
&= u^{(0,t_0)}  Z^{(\varepsilon_0)}
\left( u^{(0,t_1)} \cdots u^{({r-2},t_{r-1})}
B^{({\bm \varepsilon}')}({\bm a}', {\bm j}') ; (a_0,j_0)\right) \\
&= u^{(0,t_0)}  Z^{(\varepsilon_0)}
\left( B^{({\bm \varepsilon}')}\left( ({\bm a}', {\bm j}') ; 
(t_1, \dots, t_{r-1}) \right) ; (a_0,j_0)\right) \\
&= u^{(0,t_0)} B^{(\varepsilon_0)}(a_0,j_0) {\rm Fr'}(z) \\
&= B^{(\varepsilon_0)}\left( (a_0,j_0) ; t_0\right){\rm Fr'}(z).
\end{align*} 
Therefore, we see that  
$$B^{({\bm \varepsilon})}\left( ({\bm a}, {\bm j}); (t_0, \dots, t_{r-1}) \right) \neq 0  
\Longleftrightarrow \mbox{$B^{(\varepsilon_0)}\left( (a_0,j_0) ; t_0\right) \neq 0$ and 
${\rm Fr'}(z) \neq 0$}$$ since the multiplication map 
$\mathcal{U}_1 \otimes_k {\rm Fr}'(\mathcal{U}) \rightarrow \mathcal{U}$ 
is a $k$-linear isomorphism (see the remark of 
\cite[Proposition 2.3]{yoshii17}). From the result for $r=1$, we have 
$$B^{(\varepsilon_0)}\left( (a_0,j_0) ; t_0\right) \neq 0 \Longleftrightarrow  
\mbox{$-\widetilde{n}^{(\varepsilon_0+1)}(a_0,j_0) \leq t_0 \leq n^{(\varepsilon_0+1)}(a_0,j_0)$.}
$$ 
Moreover, by Proposition 3.4 (ii) and induction on $r$, we have 
\begin{align*}
{\rm Fr}'(z) \neq 0 & \Longleftrightarrow  
B^{({\bm \varepsilon}')}\left( ({\bm a}', {\bm j}'); (t_1, \dots, t_{r-1}) \right) \neq 0 \\ 
& \Longleftrightarrow 
\mbox{$-\widetilde{n}^{(\varepsilon_i+1)}(a_i,j_i) \leq t_i \leq n^{(\varepsilon_i+1)}(a_i,j_i)$ 
for each $i \in \{1,\dots, r-1\}.$}
\end{align*}  
Therefore, the proposition follows. $\square$ \\

For later use, we shall prove the following lemma for $\mathcal{U}_1$. \\

\begin{Lem}
Let $z \in \mathcal{U}_1$ be a nonzero element of 
degree $t$ with $-(p-1) \leq t \leq p-1$. For $(a,j) \in \mathcal{P}_{\mathbb{Z}}$,  the following hold. 
\\ \\ {\rm (i)} $z B^{(0)}(a,j)$ can be written as a $k$-linear combination of 
$B^{(0)}\left( (a,j) ; t\right)$ and $B^{(1)}\left( (a,j) ; t\right)$. \\ 
{\rm (ii)} $z B^{(1)}(a,j)$ can be written as  a 
scalar multiple  of  $B^{(1)}\left( (a,j) ; t\right)$.
\end{Lem}

\noindent {\itshape Proof.} It is enough to show (i) and (ii) when 
$$z= Y^{(n_1)} {H \choose n_2} X^{(n_3)},$$
where $0 \leq n_1,n_2,n_3 \leq p-1$ and $t=n_3-n_1$. Note that 
\begin{align*}
Y^{(n_1)} {H \choose n_2} X^{(n_3)} B^{(\varepsilon)}(a,j) 
&= Y^{(n_1)}  X^{(n_3)} {H + 2n_3 \choose n_2} B^{(\varepsilon)}(a,j) \\
&= \dfrac{1}{n_1! n_3!}{a + 2n_3 \choose n_2}Y^{n_1}  X^{n_3}  B^{(\varepsilon)}(a,j)
\end{align*}
for $\varepsilon \in \mathbb{F}_2$. Then it is easy to see that (i) and (ii) hold if $p=2$. 
So assume that $p$ is odd. 

(i) If  $j=0$, then the result is equivalent to (ii) since 
$B^{(0)}(a,0)=B^{(1)}(a,0)$. So we may assume that $j \neq 0$. Suppose that $t < 0$. By Proposition 4.3 (i) and (iii) (together with its remark) we have 
\begin{align*}
Y^{n_1}X^{n_3} B^{(0)}(a,j) 
&= Y^{n_1-n_3} Y^{n_3} X^{n_3} B^{(0)}(a,j) \\
&= Y^{n_1-n_3} \left( \beta_{n_3}(a,j) B^{(0)}(a,j)
+ 4j^2 \sum_{i=0}^{n_3-1} \dfrac{\beta_{n_3}(a,j)}{\gamma_i(a,j)} B^{(1)}(a,j)\right) \\
&= \beta_{n_3}(a,j) B^{(0)}\left( (a,j);t\right)
+ 4j^2 \sum_{i=0}^{n_3-1} \dfrac{\beta_{n_3}(a,j)}{\gamma_i(a,j)} 
B^{(1)}\left( (a,j);t \right)
\end{align*}
if $0 \leq n_3 \leq n^{(0)}(a,j)$ and 
\begin{align*}
Y^{n_1}X^{n_3} B^{(0)}(a,j) 
&= Y^{n_1-n_3} Y^{n_3} X^{n_3} B^{(0)}(a,j) \\
&= Y^{n_1-n_3} \cdot 4j^2 
\left( \prod_{i=0,\ i \neq {n}^{(0)}(a,j)}^{n_3-1}\gamma_i(a,j) \right) B^{(1)}(a,j) \\
&= 4j^2 \left( \prod_{i=0,\ i \neq {n}^{(0)}(a,j)}^{n_3-1}\gamma_i(a,j) \right) 
B^{(1)}\left( (a,j); t \right) 
\end{align*}
if $n^{(0)}(a,j) < n_3 \leq p-1$. On the other hand, suppose that 
$t \geq 0$. 
By Propositions 4.2 (i) and 4.3 (i) and (iii) (together with its remark)  we have 
\begin{align*}
\lefteqn{Y^{n_1}X^{n_3} B^{(0)}(a,j)= Y^{n_1}X^{n_1}X^t B^{(0)}(a,j) 
=Y^{n_1}X^{n_1}B^{(0)}(a+2t,j) X^t} \\
&= \left( \beta_{n_1}(a+2t,j)B^{(0)}(a+2t,j) + 4j^2 
\sum_{i=0}^{n_1-1} \dfrac{\beta_{n_1}(a+2t,j)}{\gamma_i(a+2t,j)} B^{(1)}(a+2t,j)\right) 
X^t \\
&= \beta_{n_1}(a+2t,j)B^{(0)}\left( (a,j);t \right)+ 4j^2 
\sum_{i=0}^{n_1-1} \dfrac{\beta_{n_1}(a+2t,j)}{\gamma_i(a+2t,j)} B^{(1)}\left( (a,j);t\right)
\end{align*}
if $0 \leq n_1 \leq n^{(0)}(a+2t,j)$ and 
\begin{align*}
Y^{n_1}X^{n_3} B^{(0)}(a,j)
&= Y^{n_1}X^{n_1} B^{(0)}(a+2t,j) X^t \\
&= 4j^2 \left( \prod_{i=0,\ i \neq {n}^{(0)}(a+2t,j)}^{n_1-1}\gamma_i(a+2t,j)\right) 
B^{(1)}\left( (a,j);t\right)
\end{align*}
if $n^{(0)}(a+2t,j) < n_1 \leq p-1$, as required. 

(ii) By Propositions 4.2 (i) and 4.3 (ii) (together with its remark) we have 
\begin{align*}
Y^{n_1}X^{n_3} B^{(1)}(a,j)
&= Y^{n_1-n_3} Y^{n_3} X^{n_3} B^{(1)}(a,j) \\
&= Y^{n_1-n_3} \cdot \beta_{n_3}(a,j) B^{(1)}(a,j) \\
&= \beta_{n_3}(a,j) B^{(1)}\left( (a,j); t \right)
\end{align*}
if $-(p-1) \leq t < 0$ and 
\begin{align*}
Y^{n_1}X^{n_3} B^{(1)}(a,j)
&=  Y^{n_1}X^{n_1} X^t B^{(1)}(a,j) \\
&=  Y^{n_1}X^{n_1} B^{(1)}(a+2t,j) X^t \\
&= \beta_{n_1}(a+2t,j) B^{(1)}(a+2t,j)X^t \\
&= \beta_{n_1}(a+2t,j) B^{(1)}\left( (a,j); t \right)
\end{align*}
if $0 \leq t \leq p-1$, as required. Now the lemma follows. $\square$ \\

The following lemma will be used to 
prove Theorem 5.3. \\

\begin{Lem} 
For  each $i \in \mathbb{Z}$ with $0 \leq i \leq r-1$, let $s_i$ and $t_i$ be integers satisfying 
$-(p-1) \leq s_i \leq p-1$ and  $-(p-1) \leq t_i \leq p-1$. Then, for 
$({\bm a}, {\bm j}) = \left( (a_i,j_i)\right)_{i=0}^{r-1} \in \mathcal{P}_{\mathbb{Z}}^r$, the element 
$$u^{(0,s_0)} u^{(1,s_1)} \cdots u^{({r-1}, s_{r-1})} 
B^{({\bm 1})} \left( ({\bm a}, {\bm j}); (t_0, \dots, t_{r-1})\right)$$
is a scalar multiple of 
$$B^{({\bm 1})} \left( ({\bm a}, {\bm j}); (s_0+t_0, \dots, s_{r-1}+t_{r-1})\right).$$
Moreover, if $-\widetilde{n}^{(0)}(a_i,j_i) \leq t_i \leq n^{(0)}(a_i,j_i)$ for each 
$i \in \{ 0, \dots, r-1\}$, the scalar can be taken to be  nonzero. 
\end{Lem}

\noindent {\itshape Proof.} If there exists an integer $i$ such that 
$-(p-1) \leq t_i < -\widetilde{n}^{(0)}(a_i,j_i)$ or $n^{(0)}(a_i,j_i) < t_i \leq p-1$, 
we have $B^{({\bm 1})} \left( ({\bm a}, {\bm j}); (t_0, \dots, t_{r-1})\right)=0$ by 
Proposition 4.5 and so we can take $0$ as the scalar. So we may assume that 
$-\widetilde{n}^{(0)}(a_i,j_i) \leq t_i \leq n^{(0)}(a_i,j_i)$ for each 
$i \in \{ 0, \dots, r-1\}$ (hence 
$B^{({\bm 1})} \left( ({\bm a}, {\bm j}); (t_0, \dots, t_{r-1})\right) \neq 0$).   

Consider the case of $r=1$. It is easy when $p=2$, so assume that 
$p$ is odd.  We use Propositions 4.2 (i) and 4.3 (ii) 
(together with its remark). If   $s_0 \geq 0$ and $t_0 \geq 0$ or if  
$s_0 \leq 0$ and $t_0 \leq 0$, clearly we have 
$$u^{(0,s_0)} B^{(1)}\left( (a_0,j_0) ; t_0 \right)= B^{(1)}\left( (a_0,j_0) ; s_0+t_0 \right).$$
If $s_0 \geq 0$ and $t_0 <0$, setting $t=-t_0$ we have 
\begin{align*}
u^{(0,s_0)} B^{(1)}\left( (a_0,j_0) ; t_0 \right) 
&= X^{s_0} Y^t B^{(1)}(a_0,j_0) \\
&= X^{s_0-t} X^t Y^t B^{(1)}(a_0,j_0) \\
&= \widetilde{\beta}_t(a_0,j_0) B^{(1)}\left( (a_0,j_0); s_0+t_0\right),
\end{align*}
$$\widetilde{\beta}_t(a_0,j_0) = \prod_{i=0}^{t-1} \widetilde{\gamma}_i (a_0,j_0) \neq 0$$
when $s_0 \geq t$ and 
\begin{align*}
u^{(0,s_0)} B^{(1)}\left( (a_0,j_0) ; t_0 \right) 
&= X^{s_0} Y^t B^{(1)}(a_0,j_0) \\
&= X^{s_0} Y^{s_0} B^{(1)}\left(a_0-2(t-s_0),j_0 \right) Y^{t-s_0} \\
&= \widetilde{\beta}_{s_0} (a_0-2t+2s_0,j_0)  B^{(1)}\left( (a_0,j_0); s_0+t_0\right),
\end{align*}
$$\widetilde{\beta}_{s_0} (a_0-2t+2s_0,j_0) 
= \prod_{i=0}^{s_0-1} \widetilde{\gamma}_i (a_0-2t+2s_0, j_0) 
= \prod_{i=0}^{s_0-1} \widetilde{\gamma}_{i+t-s_0} (a_0, j_0) \neq 0$$
when $s_0 < t$. On the other hand, if $s_0 <0$ and $t_0 \geq 0$, setting 
$s=-s_0$ we have 
\begin{align*}
u^{(0,s_0)} B^{(1)}\left( (a_0,j_0) ; t_0 \right) 
&= Y^s X^{t_0} B^{(1)}(a_0,j_0) \\
&= Y^{s-t_0} Y^{t_0} X^{t_0} B^{(1)}(a_0,j_0) \\
&= \beta_{t_0}(a_0,j_0) B^{(1)}\left( (a_0,j_0) ; s_0+t_0 \right),
\end{align*}
$$\beta_{t_0}(a_0,j_0) = \prod_{i=0}^{t_0-1} \gamma_i (a_0,j_0) \neq 0$$
when $s \geq t_0$ and 
\begin{align*}
u^{(0,s_0)} B^{(1)}\left( (a_0,j_0) ; t_0 \right) 
&= Y^s X^{t_0} B^{(1)}(a_0,j_0) \\
&= Y^s X^{s} B^{(1)}\left(a_0+2(t_0-s),j_0\right) X^{t_0-s} \\
&= \beta_s ( a_0+2t_0-2s,j_0) B^{(1)}\left( (a_0,j_0) ; s_0+t_0 \right),
\end{align*}
$$\beta_{s} (a_0+2t_0-2s,j_0) 
= \prod_{i=0}^{s-1} \gamma_i (a_0+2t_0-2s, j_0) 
= \prod_{i=0}^{s-1} \gamma_{i+t_0-s} (a_0, j_0) \neq 0$$
when $s< t_0$. Thus, we have shown the lemma when $r=1$.  

Let $p$ be an arbitrary prime number again and suppose that 
$r \geq 2$. By Proposition 3.4 (ii), note that 
$$Z^{(1)} \left( B^{({\bm 1}')} \left( ({\bm a}', {\bm j}'); (t_1, \dots, t_{r-1})
\right) ; (a_0,j_0)\right)=B^{(1)}(a_0,j_0) {\rm Fr}'(z)
$$
for some nonzero $z \in \mathcal{U}$, where 
${\bm 1}'=(1, \dots, 1) \in \mathbb{F}_2^{r-1}$ and 
$({\bm a}', {\bm j}') = \left( (a_i,j_i)\right)_{i=1}^{r-1} \in \mathcal{P}_{\mathbb{Z}}^{r-1}$. 
Then we have 
\begin{align*}
\lefteqn{u^{(0,s_0)} u^{(1,s_1)} \cdots u^{({r-1}, s_{r-1})} 
B^{({\bm 1})}\left( ({\bm a}, {\bm j}) ; (t_0, \dots, t_{r-1})\right)} \\
&= u^{(0,s_0)} u^{(1,s_1)} \cdots u^{({r-1}, s_{r-1})} u^{(0,t_0)}
Z^{(1)} \left( B^{({\bm 1}')} \left( ({\bm a}', {\bm j}'); (t_1, \dots, t_{r-1})
\right) ; (a_0,j_0)\right) \\
&= u^{(0,s_0)} u^{(1,s_1)} \cdots u^{({r-1}, s_{r-1})} u^{(0,t_0)} B^{(1)}(a_0,j_0) {\rm Fr}'(z).
\end{align*}
By using Proposition 4.4 and its remark repeatedly 
(if necessary), we have 
\begin{align*}
u^{(1,s_1)} \cdots u^{({r-1}, s_{r-1})} u^{(0,t_0)}B^{(1)}(a_0,j_0) 
&= u^{(1,s_1)} \cdots u^{({r-2}, s_{r-2})} u^{(0,t_0)}u^{({r-1}, s_{r-1})}B^{(1)}(a_0,j_0) \\
&= u^{(1,s_1)} \cdots u^{({r-2}, s_{r-2})} u^{(0,t_0)}B^{(1)}(a_0,j_0) u^{({r-1}, s_{r-1})} \\
&= \cdots \\
&= u^{(0,t_0)} B^{(1)}(a_0,j_0)u^{(1,s_1)} \cdots u^{({r-1}, s_{r-1})} \\
&= u^{(0,t_0)} u^{(1,s_1)} \cdots u^{({r-1}, s_{r-1})} B^{(1)}(a_0,j_0).
\end{align*}
Note that  there is a nonzero scalar $c_2 \in k$ satisfying 
\begin{align*}
\lefteqn{u^{(0,s_1)} u^{(1,s_2)} \cdots u^{({r-2},s_{r-1})} 
B^{({\bm 1}')} \left( ({\bm a}', {\bm j}'); (t_1, \dots, t_{r-1})\right)} \\ 
&= c_2 B^{({\bm 1}')} \left( ({\bm a}', {\bm j}'); (s_1+t_1, \dots, s_{r-1}+t_{r-1})\right)
\end{align*}
by induction on $r$ and that  there is a nonzero scalar $c_1 \in k$ satisfying  
$$u^{(0,s_0)} B^{(1)} \left( (a_0,j_0); t_0\right) = c_1 B^{(1)} \left( (a_0,j_0); s_0+t_0\right)$$
by the result of $r=1$. 
Then by Proposition 3.4 (ii) and (iii) we obtain 
\begin{align*}
\lefteqn{u^{(0,s_0)} u^{(1,s_1)} \cdots u^{({r-1}, s_{r-1})} 
B^{({\bm 1})}\left( ({\bm a}, {\bm j}) ; (t_0, \dots, t_{r-1})\right)} \\
&=  u^{(0,s_0)} u^{(1,s_1)} \cdots u^{({r-1}, s_{r-1})} u^{(0,t_0)} B^{(1)}(a_0,j_0) {\rm Fr}'(z) \\
&=  u^{(0,s_0)} u^{(0,t_0)}u^{(1,s_1)} \cdots u^{({r-1}, s_{r-1})}  B^{(1)}(a_0,j_0) {\rm Fr}'(z) \\
&=  u^{(0,s_0)} u^{(0,t_0)}u^{(1,s_1)} \cdots u^{({r-1}, s_{r-1})} 
Z^{(1)} \left( B^{({\bm 1}')} \left( ({\bm a}', {\bm j}'); (t_1, \dots, t_{r-1})\right) ; 
(a_0,j_0)\right)  \\
&= u^{(0,s_0)} u^{(0,t_0)}Z^{(1)} \left( u^{(0,s_1)} \cdots u^{({r-2}, s_{r-1})} 
B^{({\bm 1}')} \left( ({\bm a}', {\bm j}'); (t_1, \dots, t_{r-1})\right) ; (a_0,j_0)\right)  \\
&= u^{(0,s_0)} u^{(0,t_0)}Z^{(1)} \left( c_2
B^{({\bm 1}')} \left( ({\bm a}', {\bm j}'); 
(s_1+t_1, \dots, s_{r-1}+t_{r-1})\right) ; (a_0,j_0)\right)  \\
&= c_2 u^{(0,s_0)} u^{(0,t_0)} 
Z^{(1)} \left( B^{({\bm 1}')} \left( ({\bm a}', {\bm j}'); 
(s_1+t_1, \dots, s_{r-1}+t_{r-1})\right) ; (a_0,j_0)\right)  \\
&= c_2 u^{(0,s_0)} u^{(0,t_0)} B^{(1)}(a_0,j_0) {\rm Fr}' (z') \\
&= c_1 c_2  B^{(1)}\left( (a_0,j_0); s_0+t_0 \right) {\rm Fr}' (z') \\
&= c_1 c_2 u^{(0,s_0+t_0)}Z^{(1)} \left( 
B^{({\bm 1}')} \left( ({\bm a}', {\bm j}'); 
(s_1+t_1, \dots, s_{r-1}+t_{r-1})\right) ; (a_0,j_0)\right)  \\
&= c_1 c_2 B^{({\bm 1})} \left( ({\bm a}, {\bm j}); 
(s_0+t_0,  \dots, s_{r-1}+t_{r-1})\right) 
\end{align*} 
for some  $z' \in \mathcal{U}$ (which is not necessarily nonzero), as required. $\square$ \\

Let $V$ be a finite-dimensional $k$-subspace of $\mathcal{U}$. For  
$(a,j) \in \mathcal{P}_{\mathbb{Z}}$ and $\varepsilon \in \mathbb{F}_2$, we denote 
by  $\mathcal{U}_1 Z^{(\varepsilon)} \left( V ; (a,j)\right)$ 
the $k$-subspace spanned by all elements of the form 
$z Z^{(\varepsilon)} \left( v; (a,j)\right)$ with $z \in \mathcal{U}_1$ and 
$v \in V$.    \\

\begin{Lem}
Let $V$ be a finite-dimensional $k$-subspace of $\mathcal{U}$. Let   
$\{ v_t\}$ be  a $k$-basis of $V$.   For  
$(a,j) \in \mathcal{P}_{\mathbb{Z}}$ and $\varepsilon \in \mathbb{F}_2$, let 
$\{ w_s\}$ be a $k$-basis of   $\mathcal{U}_1 B^{(\varepsilon)}(a,j)$. Then the set 
$\left\{ w_s Z^{(0)} \left( v_t; (a,j)\right)\right\}$ forms a $k$-basis of  
$\mathcal{U}_1 Z^{(\varepsilon)} \left( V ; (a,j)\right)$. 
\end{Lem}

\noindent {\itshape Proof.} Suppose that 
$$\sum_{s,t} \alpha_{s,t} w_s Z^{(0)}\left( v_t; (a,j)\right)=0,$$
where $\alpha_{s,t} \in k$. Note that the left-hand side is equal to 
$$\sum_s w_s Z^{(0)} \left( \sum_t \alpha_{s,t} v_t ; (a,j)\right).$$
It follows from Proposition 3.4 (ii) that for each $s$, there exists 
$z_s \in \mathcal{U}$ such that 
$$Z^{(0)} \left( \sum_t \alpha_{s,t} v_t ; (a,j)\right) = 
B^{(0)}(a,j) {\rm Fr}'(z_s) $$ and hence 
$$0= \sum_s w_s Z^{(0)} \left( \sum_t \alpha_{s,t} v_t ; (a,j)\right) = 
\sum_s w_s B^{(0)}(a,j) {\rm Fr}'(z_s) = \sum_s w_s {\rm Fr}'(z_s). $$ 
Since the multiplication map 
$\mathcal{U}_1 \otimes_k {\rm Fr}'(\mathcal{U}) \rightarrow \mathcal{U}$ is a $k$-linear 
isomorphism   and since $w_s$ are linearly independent over $k$, we must have 
${\rm Fr}'(z_s)=0$ and hence 
$\sum_t \alpha_{s,t} v_t=0$ 
for each $s$ by Proposition 3.4 (i). Then the linear independence of $v_t$ implies 
$\alpha_{s,t}=0$. Therefore, the elements 
$w_s Z^{(0)}\left( v_t; (a,j)\right)$ are linearly independent over $k$. 

Finally, consider an element $zZ^{(\varepsilon)}\left( v; (a,j)\right)$ for 
$z \in \mathcal{U}_1$ and $v \in V$. If we write 
$z B^{(\varepsilon)}(a,j)= \sum_s \alpha_s w_s$ and 
$v = \sum_t \beta_{t} v_t$ with $\alpha_s, \beta_t \in k$, then the linearity of 
the map $Z^{(0)} \left( -; (a,j)\right)$ and the fact that 
$B^{(\varepsilon)}(a,j)B^{(0)}(a,j)=B^{(\varepsilon)}(a,j)$ imply 
$$zZ^{(\varepsilon)}\left( v; (a,j)\right)= 
z B^{(\varepsilon)}(a,j) Z^{(0)}\left( v; (a,j)\right)=
\sum_{s,t} \alpha_s \beta_t w_s Z^{(0)} \left( v_t; (a,j)\right).$$
Therefore, the lemma follows. $\square$

\section{Main results}

Note from Proposition 4.5 that  all  elements in  
$\widehat{\mathcal{B}}_r \left(({\bm a}, {\bm j}), {\bm \varepsilon} \right)$ are nonzero. 

The following theorem is one of the main results in this paper. 
\\

\begin{The}
Let $({\bm a}, {\bm j})  \in 
\mathcal{P}_{\mathbb{Z}}^r$ and 
${\bm \varepsilon} \in 
\mathcal{Y}_r({\bm a}, {\bm j})$. Then the set 
$\widehat{\mathcal{B}}_r \left(({\bm a}, {\bm j}), {\bm \varepsilon} \right)$ 
forms a $k$-basis of the $\mathcal{U}_r$-module 
$\mathcal{U}_r  B^{(\bm \varepsilon)}({\bm a}, {\bm j})$. 
\end{The}

Before proving the theorem, we deal with the case of $r=1$. \\

\begin{Prop}
Let $(a,j) \in \mathcal{P}_{\mathbb{Z}}$. The following hold. \\ \\
{\rm (i)}  If $(a,j)$ does not satisfy {\rm (E)}, 
the set $\widehat{\mathcal{B}}_1\left((a,j),0 \right)$ forms a $k$-basis of the 
$\mathcal{U}_1$-module $\mathcal{U}_1  B^{(0)}(a,j)$. \\ \\
{\rm (ii)} The set $\widehat{\mathcal{B}}_1\left((a,j),1 \right)$ forms a $k$-basis of the 
$\mathcal{U}_1$-module $\mathcal{U}_1  B^{(1)}(a,j)$.
\end{Prop}

\noindent {\itshape Proof.} We shall prove only (i). The argument of (ii) is similar and 
even easier. 

First let us show the $k$-linear independence of the set 
$\widehat{\mathcal{B}}_1\left((a,j),0 \right)$ in (i). 
Suppose that the sum 
$$
\sum_{\left( \theta, t(\theta) \right) \in \widehat{\Theta}_1 \left( 
(a,j), 0\right)} \alpha_{\left( \theta, t(\theta) \right)} 
B^{(\theta)}\left( (a,j);t(\theta)\right) 
= \sum_{\theta=0}^{1}\ \ \sum_{t(\theta)=-\widetilde{n}^{(\theta+1)}(a,j)}^{n^{(\theta+1)}(a,j)}
\alpha_{\left( \theta, t(\theta) \right)} B^{(\theta)}\left( (a,j);t(\theta)\right) 
$$
is equal to $0$, where $\alpha_{\left( \theta, t(\theta) \right)} \in k$. 
Since $n^{(0)}(a,j)<n^{(1)}(a,j)$ and 
$\widetilde{n}^{(0)}(a,j)<\widetilde{n}^{(1)}(a,j)$ by assumption of $(a,j)$, 
taking the part of degree $t$ we have 
$$\alpha_{(0,t)}B^{(0)}\left((a,j); t \right)+\alpha_{(1,t)}B^{(1)}\left((a,j); t \right)=0$$
for $-\widetilde{n}^{(0)}(a,j) \leq t \leq {n}^{(0)}(a,j)$ and 
$$\alpha_{(0,t)}B^{(0)}\left((a,j); t \right)=0$$
for $-\widetilde{n}^{(1)}(a,j) \leq t < -\widetilde{n}^{(0)}(a,j)$ or 
${n}^{(0)}(a,j) < t \leq {n}^{(1)}(a,j)$. In the latter case, clearly  we have 
$\alpha_{(0,t)}=0$. In the former case, since 
$\left( B^{(1)}(a,j) \right)^2=0$ by Remark (e) of Definition 4.1, 
 right multiplication by $B^{(1)}(a,j)$ implies 
$\alpha_{(0,t)}B^{(1)}\left((a,j); t \right)=0$ and hence implies $\alpha_{(0,t)}=0$ and 
$\alpha_{(1,t)}=0$. Thus, 
the claim follows. 

Let $V$ be the $k$-subspace spanned by 
$\widehat{\mathcal{B}}_1\left((a,j),0 \right)$. We wish to show that 
$V=\mathcal{U}_1  B^{(0)}(a,j)$. To show that 
$\mathcal{U}_1  B^{(0)}(a,j) \subseteq V$, it is enough to check  that 
$V$ is closed under the left action by $\mathcal{U}_1$ since $B^{(0)}(a,j) \in V$. 
For this, we have to check that 
$zB^{(\theta)}\left((a,j);t(\theta)\right) \in V$
for any $z \in \mathcal{U}_1$ and 
$\left( \theta, t(\theta) \right)  \in \widehat{\Theta}_1\left((a,j),0 \right)$.  
But without loss of generality we may assume that $z$ has 
 degree $s$ with $-(p-1) \leq s \leq p-1$. 
Then $zB^{(0)}\left((a,j);t(0)\right)$ is equal to $0$ if $zu^{(0,t(0))}=0$ and is  
a $k$-linear combination of 
$B^{(0)}\left( (a,j); s+t(0)\right)$ and  $B^{(1)}\left( (a,j); s+t(0)\right)$ if 
$zu^{(0,t(0))} \neq 0$ by Lemma 4.6 (i), and  $zB^{(1)}\left((a,j);t(1)\right)$ is equal to 
$0$ if $zu^{(0,t(1))}=0$ and is a nonzero scalar multiple of   $B^{(1)}\left( (a,j); s+t(1)\right)$ 
if $zu^{(0,t(1))} \neq 0$ by Lemma 4.6 (ii). 
Note that each 
$B^{(\rho)}\left( (a,j); t \right)$ with $\rho \in \mathbb{F}_2$ and $t \in \mathbb{Z}$ 
either lies in $\widehat{\mathcal{B}}_1\left((a,j),0 \right)$
or is  equal to zero   by Proposition 4.5. 
Thus we have 
$zB^{(\theta)}\left((a,j);t(\theta)\right) \in V$ and 
hence $\mathcal{U}_1 B^{(0)}(a,j) \subseteq V$. 

Finally, choose an element $B^{(\theta)}\left((a,j);t(\theta)\right) \in 
\widehat{\mathcal{B}}_1\left((a,j),0 \right)$ 
arbitrarily, where $\left( \theta, t(\theta) \right)  \in 
\widehat{\Theta}_1\left((a,j),0 \right)$. 
Then we have  
$$B^{(\theta)}\left( (a,j);t(\theta)\right)=u^{(0,t(\theta))}B^{(\theta)}(a,j) B^{(0)}(a,j)
\in   \mathcal{U}_1 B^{(0)}(a,j).$$  
Therefore, we obtain $V \subseteq \mathcal{U}_1  B^{(0)}(a,j)$ and (i) follows. $\square$ \\

Let us now turn to the proof of Theorem 5.1. \\

\noindent {\itshape Proof of Theorem 5.1.}  
We have already proved the case of $r=1$ in Proposition 5.2. Assume that 
$r \geq 2$. 
Write $({\bm a}, {\bm j}) = \left( (a_i,j_i)\right)_{i=0}^{r-1}$ and 
${\bm \varepsilon}= (\varepsilon_0, \dots, \varepsilon_{r-1})$. 
To apply Lemma 4.8,  we shall show that 
$$\mathcal{U}_r B^{({\bm \varepsilon})}({\bm a}, {\bm j})
= \mathcal{U}_1 Z^{(\varepsilon_0)} \left( \mathcal{U}_{r-1} B^{({\bm \varepsilon}')}
({\bm a}', {\bm j}'); (a_0,j_0)\right),  \eqno{(*)}$$
where ${\bm \varepsilon}'=(\varepsilon_1, \dots, \varepsilon_{r-1})$ and 
$({\bm a}', {\bm j}')= \left( (a_i, j_i)\right)_{i=1}^{r-1}$. 
Let $z_1 \in \mathcal{U}_1$ and $z_2 \in \mathcal{U}_{r-1}$. Since 
the algebra $\mathcal{U}_{r-1}$ is generated by $1$, $X^{(p^i)}$, and 
$Y^{(p^i)}$ with $0 \leq i \leq r-2$, there is an element $z'_2$ in the subalgebra of 
$\mathcal{U}_{r}$ generated by   $1$, $X^{(p^i)}$, and 
$Y^{(p^i)}$ with $1 \leq i \leq r-1$ such that ${\rm Fr}(z'_2)=z_2$. 
By Proposition 3.4 (iii) we have 
$$z_1 Z^{(\varepsilon_0)}\left( z_2 B^{({\bm \varepsilon}')}({\bm a}', {\bm j}') ;(a_0,j_0)\right) 
= z_1 z'_2 B^{({\bm \varepsilon})}({\bm a}, {\bm j})
\in \mathcal{U}_r B^{({\bm \varepsilon})}({\bm a}, {\bm j}).$$

On the other hand, we have to check  that $zB^{({\bm \varepsilon})}({\bm a}, {\bm j})$ lies 
in the right-hand side of $(*)$ for each $z \in \mathcal{U}_r$. Note that 
$Y^{(m_1)} X^{(m_2)} \mu_{m_3}$ with $0 \leq m_i \leq p-1$ for 
$i \in \{ 1,2,3\}$ form a $k$-basis of $\mathcal{U}_1$ 
and that $Y^{(n_1)}  \mu_{n_2}^{(r-1)}X^{(n_3)}$ with $0 \leq n_i \leq p^{r-1}-1$ for 
$i \in \{ 1,2,3\}$ form a $k$-basis of $\mathcal{U}_{r-1}$. Since the multiplication map 
$\mathcal{U}_1 \otimes_k {\rm Fr}'(\mathcal{U}_{r-1}) \rightarrow \mathcal{U}_r$ is a 
$k$-linear isomorphism (see \cite[Proposition 2.3]{yoshii17}), we may assume that 
\begin{align*}
z &=  Y^{(m_1)} X^{(m_2)} \mu_{m_3} \cdot {\rm Fr}' \left( Y^{(n_1)}  
\mu_{n_2}^{(r-1)} X^{(n_3)} \right) \\
&= Y^{(m_1)} X^{(m_2)} \mu_{m_3} \cdot Y^{(n_1p)}  
{\rm Fr}' \left( \mu_{n_2}^{(r-1)} \right)X^{(n_3p)},
\end{align*} 
where $m_i$ and $n_i$ are as above.   
By Proposition 3.3 (iv) and (v) we have 
$$
z = Y^{(m_1)} X^{(m_2)}Y^{(n_1p)}  \mu_{m_3+n_2p}^{(r)}X^{(n_3p)} 
= Y^{(m_1)} X^{(m_2)}Y^{(n_1p)}  X^{(n_3p)}\mu_{m_3+(n_2-2n_3)p}^{(r)}.
$$
Since $B^{({\bm \varepsilon})}({\bm a}, {\bm j})$ is a $\mathcal{U}_r^0$-weight vector, 
there exists $c \in k$ such that 
$$\mu_{m_3+(n_2-2n_3)p}^{(r)}B^{({\bm \varepsilon})}({\bm a}, {\bm j})
=c B^{({\bm \varepsilon})}({\bm a}, {\bm j}).$$ By Proposition 3.4 (iii) we have 
$$zB^{({\bm \varepsilon})}({\bm a}, {\bm j})= 
cY^{(m_1)} X^{(m_2)} Z^{(\varepsilon_0)} \left( Y^{(n_1)} X^{(n_3)} 
B^{({\bm \varepsilon}')}({\bm a}', {\bm j}') ;(a_0,j_0)\right),$$
which lies in $\mathcal{U}_1 Z^{(\varepsilon_0)} \left( \mathcal{U}_{r-1} B^{({\bm \varepsilon}')}
({\bm a}', {\bm j}'); (a_0,j_0)\right)$. Therefore, $(*)$ is proved. 

Now we know that the set 
$$\widehat{\mathcal{B}}_{1} \left( (a_0, j_0), \varepsilon_0 \right)
= \left\{ B^{(\theta_0)} \left( (a_0,j_0); t_0(\theta_0) \right) \ 
\left| \ \left( \theta_0, t_0(\theta_0)\right) \in 
\widehat{\Theta}_{1} \left( (a_0, j_0), \varepsilon_0\right) \right. \right\}$$
forms a basis of $\mathcal{U}_1 B^{(\varepsilon_0)}(a_0,j_0)$ by Proposition 5.2 and that 
$$\widehat{\mathcal{B}}_{r-1} \left( ({\bm a}', {\bm j}'), {\bm \varepsilon}' \right)
= \left\{ B^{({\bm \theta}')} \left( ({\bm a}',{\bm j}'); {\bm t}'({\bm \theta}') \right) 
\ \left| \ \left( {\bm \theta}', {\bm t}'({\bm \theta}')\right) \in 
\widehat{\Theta}_{r-1} \left( ({\bm a}', {\bm j}'), {\bm \varepsilon}' \right) \right. \right\}$$
forms a basis of $\mathcal{U}_{r-1} B^{({\bm \varepsilon}')}({\bm a}',{\bm j}')$ by induction 
on $r$. 
Write such ${\bm \theta}'$ and ${\bm t}'({\bm \theta}')$ as 
${\bm \theta}'= (\theta_1, \dots, \theta_{r-1})$ and 
${\bm t}'({\bm \theta}')= \left( t_1(\theta_1), \dots, t_{r-1}(\theta_{r-1})\right)$. 
By Lemma 4.8, we see that the elements 
$$
B^{(\theta_0)} \left( (a_0,j_0); t_0(\theta_0) \right) 
Z^{(0)} \left( B^{({\bm \theta}')} \left( ({\bm a}',{\bm j}'); {\bm t}'({\bm \theta}') \right) 
; (a_0,j_0)\right) $$
with $ \left( \theta_0, t_0(\theta_0)\right) \in 
\widehat{\Theta}_{1} \left( (a_0, j_0), \varepsilon_0\right)$ and 
$\left( {\bm \theta}', {\bm t}'({\bm \theta}')\right) \in 
\widehat{\Theta}_{r-1} \left( ({\bm a}', {\bm j}'), {\bm \varepsilon}' \right)$ form a basis 
of 
$$\mathcal{U}_r B^{({\bm \varepsilon})}({\bm a}, {\bm j})
= \mathcal{U}_1 Z^{(\varepsilon_0)} \left( \mathcal{U}_{r-1} B^{({\bm \varepsilon}')}
({\bm a}', {\bm j}'); (a_0,j_0)\right).$$
Set ${\bm \theta}=(\theta_0, \dots, \theta_{r-1})$ and 
${\bm t}({\bm \theta})= \left( t_0(\theta_0), \dots, t_{r-1}(\theta_{r-1})\right)$ following 
the notation in Definition 4.1 (5). 
Then we see that  
$$B^{({\bm \theta})} \left( ({\bm a},{\bm j}); {\bm t}({\bm \theta}) \right)= 
B^{(\theta_0)} \left( (a_0,j_0); t_0(\theta_0) \right) 
Z^{(0)} \left( B^{({\bm \theta}')} \left( ({\bm a}',{\bm j}'); {\bm t}'({\bm \theta}') \right) 
; (a_0,j_0)\right) $$ and that 
$ \left( \theta_0, t_0(\theta_0)\right) \in 
\widehat{\Theta}_{1} \left( (a_0, j_0), \varepsilon_0\right)$ and 
$\left( {\bm \theta}', {\bm t}'({\bm \theta}')\right) \in 
\widehat{\Theta}_{r-1} \left( ({\bm a}', {\bm j}'), {\bm \varepsilon}' \right)$ if and only if  
$\left( {\bm \theta}, {\bm t}({\bm \theta})\right) \in 
\widehat{\Theta}_{r} \left( ({\bm a}, {\bm j}), {\bm \varepsilon} \right)$. Therefore,  
the elements $B^{({\bm \theta})} \left( ({\bm a},{\bm j}); {\bm t}({\bm \theta}) \right)$ 
with $\left( {\bm \theta}, {\bm t}({\bm \theta}) \right)  \in \widehat{\Theta}_r 
\left( ({\bm a}, {\bm j}) , {\bm \varepsilon} \right)$ form a basis of 
$\mathcal{U}_r B^{({\bm \varepsilon})}({\bm a}, {\bm j})$ and 
the result follows. $\square$ \\

We would like to give some applications of Theorem 5.1.

For  $({\bm a}, {\bm j})= \left( (a_i,j_i)\right)_{i=0}^{r-1} \in \mathcal{P}_{\mathbb{Z}}^r$, define $\nu_i$ as the integer 
$$\nu_i = 
\left\{ \begin{array}{ll} 
{2j_i-1} & {\mbox{if $(a_i, j_i)$ satisfies {\rm (A)} or {\rm (D)},}} \\
{p-2j_i-1} &  {\mbox{if $(a_i, j_i)$ satisfies {\rm (B)} or {\rm (C)}}}
\end{array} \right..$$
Then we have $b_i +2n^{(0)}(a_i,j_i)=\nu_i$ for each $i \in \{ 0, \dots, r-1\}$, where 
$b_i$ are defined in Proposition 3.5 (i).

For a nonnegative integer $\lambda$, let $L(\lambda)$ be a simple $\mathcal{U}$-module 
with highest $\mathcal{U}^0$-weight $\lambda$. It is well known that  $L(\lambda)$ is also simple as a $\mathcal{U}_r$-module if 
$0 \leq \lambda \leq p^r-1$ and 
any simple $\mathcal{U}_r$-module can be obtained in this way. Each simple 
$\mathcal{U}_r$-module $L(\lambda)$ has a unique nonzero element $v_{\lambda}$ 
up to scalar multiple satisfying $X^{(n)} v_{\lambda} =0$ for any integer $n$ with 
$1 \leq n \leq p^r-1$, which is called a maximal vector. This element $v_{\lambda}$  has $\mathcal{U}^{0}_{r}$-weight 
$\lambda$: ${H \choose n} v_{\lambda} = {\lambda \choose n} v_{\lambda}$ 
for $0 \leq n \leq p^r-1$.  For details, see \cite[\S 2]{humphreys77} 
or \cite[II. \S 3]{jantzenbook}.

The following theorem says that each  $B^{({\bm 1})}({\bm a}, {\bm j})$ generates a simple $\mathcal{U}_r$-submodule of the $\mathcal{U}_r$-module $\mathcal{U}_r$.   \\

\begin{The}
Let $({\bm a}, {\bm j}) =\left( (a_i,j_i)\right)_{i=0}^{r-1} \in 
\mathcal{P}_{\mathbb{Z}}^r$. Then 
$\mathcal{U}_r  B^{({\bm 1})}({\bm a}, {\bm j})$ is a simple $\mathcal{U}_r$-module 
which is isomorphic to 
$L \left( \sum_{i=0}^{r-1} p^i \nu_i \right)$. 
\end{The}
 
\noindent {\itshape Proof.} 
Let $W$ be a simple $\mathcal{U}_r$-submodule of 
$\mathcal{U}_r  B^{({\bm 1})}({\bm a}, {\bm j})$.  Choose an arbitrary nonzero 
element $v \in W$. By Theorem 5.1 we can write $v$ as a 
$k$-linear combination of the elements 
$B^{({\bm 1})}\left( ({\bm a}, {\bm j}) ; (t_0, \dots, t_{r-1})\right)$ with 
$-\widetilde{n}^{(0)}(a_i,j_i) \leq t_i \leq {n}^{(0)}(a_i,j_i)$ for all $i \in \{0, \dots, r-1\}$. 
Choose an element 
$B^{({\bm 1})}\left( ({\bm a}, {\bm j}) ; (t_0, \dots, t_{r-1})\right)$ where its coefficient 
in the expression of $v$ is 
nonzero and $\sum_{i=0}^{r-1}p^i t_i$ is minimal. Then using Proposition 4.5 
and Lemma 4.7 we see that 
$\left( \prod_{i=0}^{r-1} X^{(p^i)n^{(0)}(a_i,j_i)-t_i} \right) v$ is a nonzero scalar multiple of 
$B^{({\bm 1})}\left( ({\bm a}, {\bm j}) ; 
\left(n^{(0)}(a_0,j_0), \dots, n^{(0)}(a_{r-1},j_{r-1})\right)\right)$ and hence that 
$$\left( \prod_{i=0}^{r-1} Y^{(p^i)n^{(0)}(a_i,j_i)}\right)
\left( \prod_{i=0}^{r-1} X^{(p^i)n^{(0)}(a_i,j_i)-t_i} \right) v$$ is a nonzero scalar multiple of 
$B^{({\bm 1})}({\bm a}, {\bm j})$. Thus  $B^{({\bm 1})}({\bm a}, {\bm j})$ lies in  $W$ 
and hence the $\mathcal{U}_r$-module $\mathcal{U}_r  B^{({\bm 1})}({\bm a}, {\bm j})$ 
must be  equal to the simple $\mathcal{U}_r$-module $W$. Moreover, since  the element 
$\left( \prod_{i=0}^{r-1} X^{(p^i)n^{(0)}(a_i,j_i)-t_i} \right) v$ 
is annihilated by the action of $X^{(p^s)}$ from the left 
for any $s \in \{0, \dots, r-1 \}$ and has $\mathcal{U}^{0}_r$-weight 
$$\sum_{i=0}^{r-1} p^i b_i + 2 \sum_{i=0}^{r-1}p^i n^{(0)}(a_i,j_i) =\sum_{i=0}^{r-1}p^i \nu_i,$$
where $b_i$ are defined in Proposition 3.5 (i),  
the $\mathcal{U}_r$-module 
$\mathcal{U}_r  B^{({\bm 1})}({\bm a}, {\bm j})$ must be 
isomorphic to $L \left( \sum_{i=0}^{r-1} p^i \nu_i \right)$. Therefore, 
the result follows. $\square$ \\ 

\noindent {\bf Remark.}  Let $Q_r(\lambda)$ be the projective cover of a simple 
$\mathcal{U}_r$-module $L(\lambda)$ with $0 \leq \lambda \leq p^r-1$. Since 
$\mathcal{U}_r   B^{({\bm 1})}({\bm a}, {\bm j})$ is a unique simple 
$\mathcal{U}_r$-submodule of 
the PIM $\mathcal{U}_r   B^{({\bm 0})}({\bm a}, {\bm j})$, 
we have 
$\mathcal{U}_r   B^{({\bm 0})}({\bm a}, {\bm j}) \cong 
Q_r\left( \sum_{i=0}^{r-1} p^i \nu_i \right)$. Therefore, we have recovered the author's result 
in \cite[Theorem 5.8]{yoshii17} without using   the $\mathcal{U}$-modules 
lifted from the PIMs $Q_r(\lambda)$.  Moreover, we also see that the sum 
$\sum_{({\bm a}, {\bm j}) \in \mathcal{P}^r}
\mathcal{U}_r   B^{({\bm 1})}({\bm a}, {\bm j})$ is the socle of the 
$\mathcal{U}_r$-module $\mathcal{U}_r$.  \\

Now we can construct a $k$-basis of the radical of the  
$\mathcal{U}_r$-module $\mathcal{U}_r B^{({\bm 0})}({\bm a}, {\bm j})$. \\

\begin{The}
Fix $({\bm a},{\bm j})= \left( (a_i,j_i)\right)_{i=0}^{r-1}\in \mathcal{P}_{\mathbb{Z}}^r$. 
Let $\mathcal{V}({\bm a}, {\bm j})$ denote the subset 
$$ 
\left\{ \left. 
B^{({\bm 0})}\left( ({\bm a}, {\bm j}); 
(t_0, \dots, t_{r-1})\right)\ \right|\ 
-\widetilde{n}^{(0)}(a_i,j_i) \leq t_i \leq n^{(0)}(a_i,j_i),\ \forall i \right\}$$
of  $\mathcal{B}_r \left(({\bm a}, {\bm j}), {\bm 0} \right)$. 
Then  its complement 
$\mathcal{B}_r \left(({\bm a}, {\bm j}), {\bm 0} \right) 
\backslash \mathcal{V}({\bm a}, {\bm j})$ 
forms a $k$-basis of the radical of the projective indecomposable $\mathcal{U}_r$-module 
$\mathcal{U}_r B^{({\bm 0})}({\bm a}, {\bm j})$, and 
the image of $\mathcal{V}({\bm a}, {\bm j})$ in  the quotient space  
$$\mathcal{U}_r B^{({\bm 0})}({\bm a}, {\bm j})/{\rm rad}_{\mathcal{U}_r}
\left( \mathcal{U}_r B^{({\bm 0})}({\bm a}, {\bm j})\right)$$  forms its $k$-basis. 
\end{The}

\noindent {\itshape Proof.} We should keep in mind that 
$B^{({\bm 0})}({\bm a}, {\bm j})=B^{({\bm \tau}({\bm a}, {\bm j}))}({\bm a}, {\bm j})$ 
and 
$\mathcal{B}_r \left(({\bm a}, {\bm j}), {\bm 0} \right)=
\widehat{\mathcal{B}}_r \left(({\bm a}, {\bm j}), {\bm \tau}({\bm a},{\bm j})\right)$.  
Consider the surjective  $\mathcal{U}_r$-module homomorphism 
$$f_{({\bm a}, {\bm j})}: \mathcal{U}_r B^{({\bm 0})}({\bm a}, {\bm j}) \rightarrow 
\mathcal{U}_r B^{({\bm 1})}({\bm a}, {\bm j})$$
defined by right multiplication by $B^{({\bm 1})} ({\bm a}, {\bm j})=
B^{({\bm \sigma}({\bm a}, {\bm j}))} ({\bm a}, {\bm j})$. 
It follows from the simplicity of the 
$\mathcal{U}_r$-module $\mathcal{U}_r B^{({\bm 1})}({\bm a}, {\bm j})$ that  
${\rm Ker} f_{({\bm a}, {\bm j})} = {\rm rad}_{\mathcal{U}_r}
\left( \mathcal{U}_r B^{({\bm 0})}({\bm a}, {\bm j})\right)$. 

Following Theorem 5.1, we write an element 
$z \in \mathcal{U}_r B^{({\bm 0})}({\bm a}, {\bm j})
=\mathcal{U}_r B^{({\bm \tau}({\bm a}, {\bm j}))}({\bm a}, {\bm j})$ as a $k$-linear 
combination of the  elements in 
$\mathcal{B}_r \left(({\bm a}, {\bm j}), {\bm 0}\right)=
\widehat{\mathcal{B}}_r \left(({\bm a}, {\bm j}), {\bm \tau}({\bm a},{\bm j})\right)$: 
$$z= \sum_{({\bm \theta},{\bm t}({\bm \theta})) \in \widehat{\Theta}_{r}
\left(({\bm a}, {\bm j}), {\bm \tau}({\bm a},{\bm j}) \right)} 
 \alpha_{({\bm \theta},{\bm t}({\bm \theta}))}
B^{({\bm \theta})}\left( ({\bm a}, {\bm j}); {\bm t}({\bm \theta}) \right) \eqno{(**)}$$ 
with $\alpha_{({\bm \theta},{\bm t}({\bm \theta}))} \in k$. Note that 
$$B^{({\bm \theta})} ({\bm a}, {\bm j})B^{({\bm \sigma}({\bm a}, {\bm j}))} 
({\bm a}, {\bm j}) = 
\left\{ 
\begin{array}{ll}
B^{({\bm 1})} ({\bm a}, {\bm j}) & \mbox{if ${\bm \theta}={\bm \tau}({\bm a},{\bm j}) $,} \\
0 & \mbox{if ${\bm \theta} \neq {\bm \tau}({\bm a},{\bm j}) $}
\end{array} 
\right.$$
for ${\bm \theta} \in \mathcal{Y}_r({\bm a},{\bm j})$ 
by Remark (e) of Definition 4.1 
(since ${\bm \sigma}({\bm a}, {\bm j}) \in \mathcal{X}_r({\bm a},{\bm j})$) and that 
$$B^{({\bm 1})}\left( ({\bm a}, {\bm j}); (t_0, \dots, t_{r-1})\right) \neq 0 
\Longleftrightarrow \mbox{$-\widetilde{n}^{(0)}(a_i,j_i) \leq t_i \leq n^{(0)}(a_i,j_i)$,\ $\forall 
i$}$$ by Proposition 4.5.
By multiplying the both sides of the equality $(**)$ by 
$B^{({\bm \sigma}({\bm a}, {\bm j}))}({\bm a}, {\bm j})$ 
from the right, we obtain  
$$f_{({\bm a}, {\bm j})}(z)= \sum_{{\bm t}({\bm \tau}({\bm a}, {\bm j}))} 
\alpha_{({\bm \tau}({\bm a}, {\bm j}),{\bm t}({\bm \tau}({\bm a}, {\bm j})) )}
B^{({\bm 1})}\left( ({\bm a}, {\bm j});{\bm t}({\bm \tau}({\bm a}, {\bm j}))  \right),$$
where ${\bm t}({\bm \tau}({\bm a}, {\bm j}))=
(t_0(\tau_0), \dots, t_{r-1}(\tau_{r-1}))$ runs through the elements 
in $\mathbb{Z}^r$ with each $t_i(\tau_i)$ satisfying  
$$-\widetilde{n}^{(0)}(a_i,j_i) \leq t_i(\tau_i) \leq n^{(0)}(a_i,j_i)$$
(recall that ${\bm \tau}({\bm a},{\bm j})=(\tau_0, \dots, \tau_{r-1})$). 
Hence we see that $f_{({\bm a}, {\bm j})}(z)=0$ if and only if 
$\alpha_{({\bm \tau}({\bm a},{\bm j}), {\bm t}({\bm \tau}({\bm a},{\bm j})))} =0$ 
in the equality $(**)$ whenever 
$$-\widetilde{n}^{(0)}(a_i,j_i) \leq t_i(\tau_i) \leq n^{(0)}(a_i,j_i)$$
for any $i$.  
Therefore, we see that the subset 
$\widehat{\mathcal{B}}_r \left(({\bm a}, {\bm j}), {\bm \tau}({\bm a},{\bm j}) \right) 
\backslash \mathcal{V}({\bm a}, {\bm j})$ 
forms a $k$-basis of ${\rm Ker} f_{({\bm a}, {\bm j})}= {\rm rad}_{\mathcal{U}_r}
\left( \mathcal{U}_r B^{({\bm 0})}({\bm a}, {\bm j})\right)$  
and that the image of $\mathcal{V}({\bm a}, {\bm j})$ under the map 
$f_{({\bm a}, {\bm j})}$ forms a $k$-basis of $\mathcal{U}_r B^{({\bm 1})}({\bm a}, {\bm j})$. 
Now the theorem follows. $\square$\\ 

Theorem 5.4  enables us to construct a $k$-basis of the Jacobson radical 
${\rm rad} \mathcal{U}_r$ and to know that each $\mathcal{U}_r$-module 
$\mathcal{U}_r B^{({\bm \varepsilon})}({\bm a}, {\bm j})$ with 
$({\bm a},{\bm j}) \in \mathcal{P}_{\mathbb{Z}}^r$ and 
${\bm \varepsilon} \in \mathbb{F}_2^{r}$ 
has simple head and 
simple socle.  \\

\begin{Cor}
Let $\mathcal{V} $ denote the subset 
$$
\left\{  
B^{({\bm 0})}\left( ({\bm a}, {\bm j}); 
(t_0, \dots, t_{r-1})\right)\ \left|\ 
\begin{array}{l}
{({\bm a},{\bm j})= \left( (a_i,j_i)\right)_{i=0}^{r-1}\in \mathcal{P}^r,} \\
{-\widetilde{n}^{(0)}(a_i,j_i) \leq t_i \leq n^{(0)}(a_i,j_i),\ \forall i} 
\end{array}
\right. \right\}$$
of  $\bigcup_{({\bm a},{\bm j}) \in \mathcal{P}^r}
\mathcal{B}_r \left(({\bm a}, {\bm j}), {\bm 0}\right)$. Then  
its complement 
$\bigcup_{({\bm a},{\bm j}) \in \mathcal{P}^r}
\mathcal{B}_r \left(({\bm a}, {\bm j}), {\bm 0}\right)
\backslash \mathcal{V}$ 
forms a $k$-basis of  ${\rm rad}\mathcal{U}_r$, and 
the image of $\mathcal{V}$ in  the quotient space  
$\mathcal{U}_r/{\rm rad}\mathcal{U}_r$  forms its $k$-basis. 
\end{Cor}

\noindent {\itshape Proof.} 
Clearly $\mathcal{V}$ is a disjoint union of all $\mathcal{V}({\bm a},{\bm j})$ with 
$({\bm a},{\bm j}) \in \mathcal{P}^r$, where $\mathcal{V}({\bm a},{\bm j})$ is the set 
defined in Theorem 5.4. 
Consider the direct sum decomposition of the 
$\mathcal{U}_r$-module $\mathcal{U}_r$ into 
PIMs induced by 
$\sum_{({\bm a},{\bm j}) \in \mathcal{P}^r} B^{({\bm 0})}({\bm a},{\bm j})=1$: 
$$\mathcal{U}_r = \sum_{({\bm a},{\bm j}) \in \mathcal{P}^r} \mathcal{U}_r 
B^{({\bm 0})}({\bm a},{\bm j})
=\sum_{({\bm a},{\bm j}) \in \mathcal{P}^r} \mathcal{U}_r 
B^{({\bm \tau}({\bm a},{\bm j}))}({\bm a},{\bm j}).$$
Thus we have 
$${\rm rad} \mathcal{U}_r = \sum_{({\bm a},{\bm j}) \in \mathcal{P}^r} 
{\rm rad}_{\mathcal{U}_r} \big( \mathcal{U}_r 
B^{({\bm \tau}({\bm a},{\bm j}))}({\bm a},{\bm j})\big).$$
Recall that $\mathcal{B}_r \left(({\bm a}, {\bm j}), {\bm 0}\right)=
\widehat{\mathcal{B}}_r \left(({\bm a}, {\bm j}), {\bm \tau}({\bm a},{\bm j})\right) $. 
It follows from Theorem 5.4 that  the subset  
$\widehat{\mathcal{B}}_r \left(({\bm a}, {\bm j}), {\bm \tau}({\bm a},{\bm j})\right) 
\backslash \mathcal{V}({\bm a}, {\bm j})$ forms a $k$-basis of the radical 
${\rm rad}_{\mathcal{U}_r} \left( \mathcal{U}_r 
B^{({\bm \tau}({\bm a},{\bm j}))}({\bm a},{\bm j})\right)$ of the 
$\mathcal{U}_r$-module $\mathcal{U}_r 
B^{({\bm \tau}({\bm a},{\bm j}))}({\bm a},{\bm j})$. Therefore, the union 
$$\bigcup_{({\bm a},{\bm j}) \in \mathcal{P}^r} \left( 
\widehat{\mathcal{B}}_r \left(({\bm a}, {\bm j}), {\bm \tau}({\bm a},{\bm j})\right) \backslash  
\mathcal{V}({\bm a}, {\bm j}) \right)=\bigcup_{({\bm a},{\bm j}) \in \mathcal{P}^r} 
 \widehat{\mathcal{B}}_r \left(({\bm a}, {\bm j}), {\bm \tau}({\bm a},{\bm j})\right) \backslash  
\mathcal{V}$$
forms a $k$-basis of ${\rm rad} \mathcal{U}_r$, and the first claim follows. 

Recall that $${\rm soc}_{\mathcal{U}_r} \mathcal{U}_r 
= \sum_{({\bm a},{\bm j}) \in \mathcal{P}^r}\mathcal{U}_r 
B^{({\bm 1})}({\bm a},{\bm j}).$$
Hence we see that there is a natural surjective $\mathcal{U}_r$-module homomorphism 
$$f= \sum_{({\bm a},{\bm j}) \in \mathcal{P}^r}f_{({\bm a},{\bm j})}: \mathcal{U}_r \rightarrow 
{\rm soc}_{\mathcal{U}_r} \mathcal{U}_r$$
induced by the maps $f_{({\bm a},{\bm j})}$ with $({\bm a},{\bm j}) \in \mathcal{P}^r$ 
defined in the proof of Theorem 5.4 and 
that ${\rm Ker} f = {\rm rad} \mathcal{U}_r $. The image of $\mathcal{V}$ under the map 
$f$ is 
$$\left\{  
B^{({\bm 1})}\left( ({\bm a}, {\bm j}); (t_0, \dots, t_{r-1})\right)\ \left|\ 
\begin{array}{l}
{({\bm a},{\bm j})= \left( (a_i,j_i)\right)_{i=0}^{r-1}\in \mathcal{P}^r,} \\
{-\widetilde{n}^{(0)}(a_i,j_i) \leq t_i \leq n^{(0)}(a_i,j_i),\ \forall i} 
\end{array}
\right. \right\},$$
which is a $k$-basis of 
$\sum_{({\bm a},{\bm j}) \in \mathcal{P}^r}\mathcal{U}_r 
B^{({\bm 1})}({\bm a},{\bm j})={\rm soc}_{\mathcal{U}_r} \mathcal{U}_r$. 
Therefore, the second claim 
follows. $\square$ \\

\noindent {\bf Remark.} Wong gives some good generating sets of $\mathcal{U}_1$ 
in \cite{wong83} if $p$ is odd.  The basis obtained in Corollary 5.5 is expected to be useful 
in improving his result. We have already had an outlook on the improvement  and 
will report it in a future paper.  \\
\

\begin{Cor}
Let $({\bm a},{\bm j}) \in \mathcal{P}_{\mathbb{Z}}^r$ and 
${\bm \varepsilon} \in \mathbb{F}_2^r$. Then the $\mathcal{U}_r$-module 
$\mathcal{U}_r B^{({\bm \varepsilon})}({\bm a}, {\bm j})$ has  head and  socle 
isomorphic to the simple $\mathcal{U}_r$-module 
$\mathcal{U}_r B^{({\bm 1})}({\bm a}, {\bm j})$. 
\end{Cor}

\noindent {\itshape Proof.} 
Without loss of generality, we may assume that  
$({\bm a},{\bm j}) \in \mathcal{P}^r$ and 
${\bm \varepsilon} \in \mathcal{Y}_r({\bm a}, {\bm j})$ (see Remark (b) of Definition 4.1). 

Since $\mathcal{U}_r B^{({\bm \varepsilon})}({\bm a}, {\bm j})$ is a $\mathcal{U}_r$-submodule 
of the PIM $\mathcal{U}_r B^{({\bm 0})}({\bm a}, {\bm j})$, clearly we have 
$${\rm soc}_{\mathcal{U}_r} \left( 
\mathcal{U}_r B^{({\bm \varepsilon})}({\bm a}, {\bm j})\right) = 
\mathcal{U}_r B^{({\bm 1})}({\bm a}, {\bm j}).$$

Consider the surjective $\mathcal{U}_r$-module homomorphism 
$$f_{({\bm a}, {\bm j})}^{({\bm \varepsilon})}: 
\mathcal{U}_r B^{({\bm \varepsilon})}({\bm a}, {\bm j}) \rightarrow 
\mathcal{U}_r B^{({\bm 1})}({\bm a}, {\bm j})$$
defined by right multiplication by $B^{({\bm 1}-{\bm \varepsilon})} ({\bm a}, {\bm j})$. 
Write $({\bm a}, {\bm j})= \left( (a_i,j_i)\right)_{i=0}^{r-1}$ and 
${\bm \varepsilon}=(\varepsilon_0, \dots, \varepsilon_{r-1})$. 
Recall from Theorem 5.1 that 
the set 
$\widehat{\mathcal{B}}_r \left( ({\bm a}, {\bm j}), {\bm \varepsilon}\right)$
forms a $k$-basis of the $\mathcal{U}_r$-module 
$\mathcal{U}_r B^{({\bm \varepsilon})}({\bm a}, {\bm j})$.

Write an element $z \in \mathcal{U}_r B^{({\bm \varepsilon})}({\bm a}, {\bm j})$ 
as a $k$-linear 
combination of the basis: 
$$z= \sum_{({\bm \theta}, {\bm t}({\bm \theta})) \in 
\widehat{\Theta}_r(({\bm a}, {\bm j}), {\bm \varepsilon})} 
\alpha_{({\bm \theta}, {\bm t}({\bm \theta}))}
B^{({\bm \theta})}\left( ({\bm a}, {\bm j}); {\bm t}({\bm \theta}) \right) \eqno{(***)}$$ 
with $\alpha_{({\bm \theta}, {\bm t}({\bm \theta}))} \in k$. Note that 
$$B^{({\bm \theta})} ({\bm a}, {\bm j})B^{({\bm 1}-{\bm \varepsilon})} ({\bm a}, {\bm j}) = 
\left\{ 
\begin{array}{ll}
B^{({\bm 1})} ({\bm a}, {\bm j}) & \mbox{if ${\bm \theta}={\bm \varepsilon}$,} \\
0 & \mbox{if ${\bm \theta}>{\bm \varepsilon}$}
\end{array} 
\right.$$
for $ {\bm \theta} \in \mathcal{Y}_r ({\bm a}, {\bm j})$ with 
${\bm \theta} \geq {\bm \varepsilon}$ 
by Remark (e) of Definition 4.1  
(since ${\bm 1}-{\bm \varepsilon} \in \mathcal{X}_r ({\bm a}, {\bm j})$) 
and that $$B^{({\bm 1})}\left( ({\bm a}, {\bm j}); (t_0,\dots, t_{r-1})\right) \neq 0 
\Longleftrightarrow \mbox{$-\widetilde{n}^{(0)}(a_i,j_i) \leq t_i \leq n^{(0)}(a_i,j_i)$, 
 $\forall i$}$$ by Proposition 4.5.
By multiplying the both sides of the equality $(***)$ by 
$B^{({\bm 1}-{\bm \varepsilon})}({\bm a}, {\bm j})$ 
from the right, we have 
$$f^{({\bm \varepsilon})}_{({\bm a}, {\bm j})}(z)= 
\sum_{{\bm t}({\bm \varepsilon})} \alpha_{({\bm \varepsilon},{\bm t}({\bm \varepsilon}))}
B^{({\bm 1})}\left( ({\bm a}, {\bm j}); {\bm t}({\bm \varepsilon}) \right),$$
where ${\bm t}({\bm \varepsilon})=(t_0(\varepsilon_0), \dots, t_{r-1}(\varepsilon_{r-1}))$ 
runs through the elements in 
$\mathbb{Z}^r$ with each $t_i(\varepsilon_i)$ satisfying  
$$-\widetilde{n}^{(0)}(a_i,j_i) \leq t_i(\varepsilon_i) \leq n^{(0)}(a_i,j_i).$$
Hence we see that $f^{({\bm \varepsilon})}_{({\bm a}, {\bm j})}(z)=0$ if and only if 
$\alpha_{({\bm \varepsilon},{\bm t}({\bm \varepsilon}))} =0$ in the equality $(***)$ whenever 
$$-\widetilde{n}^{(0)}(a_i,j_i) \leq t_i(\varepsilon_i) \leq n^{(0)}(a_i,j_i)$$
for any $i$. Therefore, we see that the subset 
$$\widehat{\mathcal{B}}_r \left(({\bm a}, {\bm j}), {\bm \varepsilon} \right) \backslash 
\left\{ \left. B^{({\bm \varepsilon})}\left( ({\bm a}, {\bm j}); 
{\bm t}({\bm \varepsilon})\right)\ \right|\ -\widetilde{n}^{(0)}(a_i,j_i) 
\leq t_i(\varepsilon_i) \leq n^{(0)}(a_i,j_i),\ \forall i \right\}$$ 
forms a $k$-basis of ${\rm Ker} f^{({\bm \varepsilon})}_{({\bm a}, {\bm j})}$. 
On the other hand, note that 
$${\rm rad}_{\mathcal{U}_r} \left( \mathcal{U}_r B^{({\bm \varepsilon})}({\bm a}, {\bm j})\right)
= ({\rm rad} \mathcal{U}_r) B^{({\bm \varepsilon})}({\bm a}, {\bm j})$$
and from Corollary 5.5 that the subset 
$\bigcup_{\left(\widetilde{\bm a}, \widetilde{\bm j}\right) \in \mathcal{P}^r}
\widehat{\mathcal{B}}_r \left(\left(\widetilde{\bm a}, \widetilde{\bm j}\right), 
{\bm \tau}\left(\widetilde{\bm a}, \widetilde{\bm j}\right)\right) \backslash \mathcal{V}$ 
forms a $k$-basis of ${\rm rad} \mathcal{U}_r$. By multiplying the elements of the subset 
by $B^{({\bm \varepsilon})}({\bm a}, {\bm j})$ from the right we see that the subset 
$$\widehat{\mathcal{B}}_r \left(({\bm a}, {\bm j}), {\bm \varepsilon} \right) \backslash 
\left\{ \left. B^{({\bm \varepsilon})}\left( ({\bm a}, {\bm j}); {\bm t}\right)\ \right|\ 
-\widetilde{n}^{(0)}(a_i,j_i) \leq t_i \leq n^{(0)}(a_i,j_i),\ \forall i \right\},$$ 
where ${\bm t}=(t_0, \dots, t_{r-1})$,  forms a $k$-basis of 
$ ({\rm rad} \mathcal{U}_r) B^{({\bm \varepsilon})}({\bm a}, {\bm j})$. 
Therefore, we have 
$${\rm Ker} f^{({\bm \varepsilon})}_{({\bm a}, {\bm j})}= 
({\rm rad} \mathcal{U}_r) B^{({\bm \varepsilon})}({\bm a}, {\bm j})=
{\rm rad}_{\mathcal{U}_r} \left( \mathcal{U}_r B^{({\bm \varepsilon})}({\bm a}, {\bm j})\right)$$
and the corollary follows. $\square$ \\ 

\noindent {\bf Remark.} Actually, all the $\mathcal{U}_r$-modules 
$\mathcal{U}_r B^{({\bm \varepsilon})}({\bm a}, {\bm j})$ are self-dual. 
This fact will be clear in a future paper. \\ \\

\noindent {\bf \Large Acknowledgment} \\

This work was supported by JSPS KAKENHI Grant Number JP18K03203.

\end{document}